\documentclass[12pt]{amsart}
\usepackage{amsmath,amssymb,latexsym,esint,cite,mathrsfs}
\usepackage{verbatim,wasysym}
\usepackage[left=2.6cm,right=2.6cm,top=3cm,bottom=3cm]{geometry}
\makeatletter \@addtoreset{equation}{section} \makeatother

\usepackage{graphicx}
\usepackage{subfigure}
\usepackage{graphicx,epic,eepic}
\usepackage{tikz}  
\usepackage{color,enumitem,graphicx}
\usepackage[colorlinks=true,urlcolor=blue,
citecolor=red,linkcolor=blue,linktocpage,pdfpagelabels,
bookmarksnumbered,bookmarksopen]{hyperref}

\numberwithin{equation}{section}
\newtheorem{theorem}{Theorem}[section]
\newtheorem{lemma}[theorem]{Lemma}

\newtheorem{proposition}[theorem]{Proposition}
\newtheorem{remark}[theorem]{Remark}
\newtheorem{corollary}[theorem]{Corollary}
\numberwithin{equation}{section}

\allowdisplaybreaks

\begin{document}

\title[Stability of (CKN) inequality]
{Gradient stability of Caffarelli-Kohn-Nirenberg inequality involving weighted $p$-Laplace}

\author[S. Deng]{Shengbing Deng$^{\ast}$}
\address{\noindent Shengbing Deng
\newline
School of Mathematics and Statistics, Southwest University,
Chongqing 400715, People's Republic of China}\email{shbdeng@swu.edu.cn}

\author[X. Tian]{Xingliang Tian}
\address{\noindent Xingliang Tian  \newline
School of Mathematics and Statistics, Southwest University,
Chongqing 400715, People's Republic of China.}\email{xltian@email.swu.edu.cn}

\thanks{$^{\ast}$ Corresponding author}

\thanks{2020 {\em{Mathematics Subject Classification.}} 35P30, 35J92, 26D10.}

\thanks{{\em{Key words and phrases.}} Caffarelli-Kohn-Nirenberg inequality; Weighted $p$-Laplace equation; Non-degeneracy; Gradient stability; Remainder term}

\allowdisplaybreaks

\begin{abstract}
The best constant and extremal functions are well known of the following Caffarelli-Kohn-Nirenberg inequality
\[
\int_{\mathbb{R}^N}|\nabla u|^p\frac{\mathrm{d}x}{|x|^{\mu}}\geq \mathcal{S}
    \left(\int_{\mathbb{R}^N}|u|^r\frac{\mathrm{d}x}{|x|^s}
    \right)^{\frac{p}{r}},
    \quad \mbox{for all}\quad  u\in C^\infty_c(\mathbb{R}^N),
\]
where $1<p<p+\mu<N$, $\frac{\mu}{p}\leq \frac{s}{r}<\frac{\mu}{p}+1$, $r=\frac{p(N-s)}{N-p-\mu}$.
An important task is investigating the stability of extremals for this inequality.
Firstly, we give the classification to the linearized problem related to the extremals which shows the extremals are non-degenerate. Then we investigate the gradient type remainder term of previous inequality by using spectral estimate combined with a compactness argument which partially extends the work of Wei and Wu [Math. Ann., 2022] to a general $p$-Laplace case, and also the work of Figalli and Zhang [Duke Math. J., 2022] to a weighted case.
\end{abstract}

\vspace{3mm}

\maketitle

\section{\bfseries Introduction}\label{sectir}

\subsection{Motivation}\label{subsectmot}
    Given $N\geq 2$ and $p\in (1,N)$, denote the homogeneous Sobolev space $\mathcal{D}^{1,p}_0(\mathbb{R}^N)$ be the closure of $C^\infty_c(\mathbb{R}^N)$ with respect to the norm
    $
    \|u\|_{\mathcal{D}^{1,p}_0(\mathbb{R}^N)}
    :=\left(\int_{\mathbb{R}^N}|\nabla u|^p\mathrm{d}x\right)^\frac{1}{p}.
    $
    The Sobolev inequality states as
    \begin{equation}\label{bzsi}
    \|\nabla u\|^p_{L^p(\mathbb{R}^N)}\geq \mathcal{S}_0\|u\|^p_{L^{p^*}(\mathbb{R}^N)},\quad \mbox{for all}\quad u\in \mathcal{D}^{1,p}_0(\mathbb{R}^N).
    \end{equation}
    for some constant $\mathcal{S}_0=\mathcal{S}_0(N,p)>0$ depending on $N$ and $p$, where $p^*:=\frac{pN}{N-p}$.
    It is well known that Aubin \cite{Au76} and Talenti \cite{Ta76} found the optimal constant and the extremals for inequality \eqref{bzsi}. Indeed, equality is achieved precisely by the functions (up to scalar multiplications)
    \begin{equation}\label{defvlz}
    V_{\lambda,z}(x):=\lambda^{\frac{N-p}{p}}V(\lambda(x-z)),\quad \mbox{for all}\quad \lambda>0,\quad z\in \mathbb{R}^N,
    \end{equation}
    where
    \begin{equation*}
    V(x)=\gamma_{N,p}(1+|x|^{\frac{p}{p-1}})^{-\frac{N-p}{p}},\quad \mbox{for some constant}\quad  \gamma_{N,p}>0. 
    \end{equation*}

    Recently, there has been a growing interest in understanding quantitative stability of inequality (\ref{bzsi}). For $p=2$, in \cite{BrL85}, Brezis and Lieb asked the question whether a remainder term - proportional to the quadratic distance of the function $u$ to be the manifold  $\mathcal{M}_0:=\{cV_{\lambda,z}: c\in\mathbb{R}, \lambda>0, z\in\mathbb{R}^N\}$ - can be added to the right hand side of (\ref{bzsi}). This question was answered affirmatively by Bianchi and Egnell \cite{BE91} by using spectral estimate combined with Lions' concentration and compactness principle (see \cite{Li85-2}), which reads that there exists constant $c_{\mathrm{BE}}>0$ such that for all $u\in \mathcal{D}^{1,2}_0(\mathbb{R}^N)$,
    \begin{equation}\label{defcbe}
    \int_{\mathbb{R}^N}|\nabla u|^2 \mathrm{d}x- \mathcal{S}_0\left(\int_{\mathbb{R}^N}|u|^{2^*} \mathrm{d}x\right)^{\frac{2}{2^*}}\geq c_{\mathrm{BE}} \inf_{v\in \mathcal{M}_0}\|\nabla u-\nabla v\|^2_{L^2(\mathbb{R}^N)}.
    \end{equation}
    After that, this result was extended later to the second-order case by Lu and Wei \cite{LW00}, arbitrary high order by Bartsch et al. \cite{BWW03}, and the whole fractional order case was proved in \cite{CFW13}. See also \cite{Ng19} for Gagliardo-Nirenberg-Sobolev inequality in quantitative form.

    While for $p\neq 2$, it needs much delicate analysis to deal with the stability of inequality (\ref{bzsi}). Cianchi et al. \cite{CFMP09} first proved a stability version of Lebesgue-type for every $1<p<N$, Figalli and Neumayer \cite{FN19} proved the gradient stability for the Sobolev inequality when $p\geq 2$, Neumayer \cite{Ne20} extended the result in \cite{FN19} to all $1<p<N$. It is worth to mention that very recently, Figalli and Zhang \cite{FZ22} obtained the sharp stability of extremals for Sobolev inequality (\ref{bzsi}) for any $1<p<N$ which reads that for all $u\in \mathcal{D}^{1,p}_0(\mathbb{R}^N)$,
    \[
    \frac{\|\nabla u\|_{L^p(\mathbb{R}^N)}}{\|u\|_{L^{p^*}(\mathbb{R}^N)}}
    -\mathcal{S}_0^{\frac{1}{p}}
    \geq c_{\mathrm{FZ}} \inf_{v\in \mathcal{M}_0}\left(\frac{\|\nabla u-\nabla v\|_{L^p(\mathbb{R}^N)}}{\|\nabla u\|_{L^p(\mathbb{R}^N)}}\right)^{\max\{2,p\}},
    \]
    for some constant $c_{\mathrm{FZ}}>0$, 
     furthermore the exponent $\max\{2,p\}$ is sharp.

    Now, let us mention the famous  Caffarelli-Kohn-Nirenberg (we write (CKN) for short) inequality which was first introduced in 1984 by Caffarelli, Kohn and Nirenberg in their celebrated work \cite{CKN84}.
    It is worth noting that many well-known and important inequalities such as Gagliardo-Nirenberg inequalities, Sobolev inequalities, Hardy-Sobolev inequalities, Nash inequalities, etc. are just special cases of (CKN) inequalities.
    Here, we just focus on the case without interpolation term, that is,
    \begin{equation}\label{cknwit}
    \left(\int_{\mathbb{R}^N}|x|^{-br}|u|^r \mathrm{d}x\right)^{\frac{p}{r}}
    \leq C_{\mathrm{CKN}}\int_{\mathbb{R}^N}|x|^{-ap}|\nabla u|^p \mathrm{d}x, \quad \mbox{for all}\quad u\in C^\infty_c(\mathbb{R}^N),
    \end{equation}
    for some constant $C_{\mathrm{CKN}}>0$, where
    \begin{equation}\label{cknwitc}
    1<p<N,\quad -\infty<a<a_c:=\frac{N-p}{p},\quad a\leq b\leq a+1,\quad r=\frac{pN}{N-p(1+a-b)}.
    \end{equation}
    A natural question is whether the best constant $C_{\mathrm{CKN}}$ could be achieved or not? Moreover, if $C_{\mathrm{CKN}}$ is achieved, are the extremal functions radial symmetry? In fact, the weights $|x|^{-ap}$ and $|x|^{-br}$ have deep influences in many aspects about this inequality, for example, achievable of best constant and symmetry of minimizers.

\begin{figure}[ht]
\begin{tikzpicture}[scale=4]
		\draw[->,ultra thick](-2.8,0)--(0,0)node[below right]{$O$}--(0.6,0)node[below]{$a$};
		\draw[->,ultra thick](0,-0.75)--(0,1.3)node[left]{$b$};
        \draw[fill=gray,domain=0:-2.8]plot(\x,{4*(1-0.47*\x)/
        (2*((1-0.47*\x)^2+3)^0.5)
        +0.47*\x -1})--(-2.8,-0.75)--(-1.5,-0.75);
        \draw[fill=brown,domain=0:-2.8]plot(\x,{4*(1-0.47*\x)/
        (2*((1-0.47*\x)^2+3)^0.5)
        +0.47*\x -1})--(-2.8,-0.5)--(0.5,1.25)--(0.5,0.25)--(0,0);
        \draw[densely dashed](0.5,0)node[below]{$a_c$}--(0.5,1.25);
        \draw[densely dashed](-2.8,-0.5)--(0.5,1.25);
        \draw[densely dashed](-1.5,-0.75)--(0,0);
        \draw[-,ultra thick](0,-0.75)--(0,1);
        \draw[-,ultra thick](-2.8,0)--(0.6,0);
        \node[left] at(0,1.05){$1$};
		\node[left] at(-0.5,0.75){$b=a+1$};
        \node[right] at (-0.4,-0.22){$b=a$};
        \node[right] at (-1.5,-0.35){$b=b_{\mathrm{FS}}(a)$};
        \node[right] at (-2.7,-0.69){$symmetry\ breaking\ region$};
        \node[right] at (-1,0.3){$symmetry\ region$};
\end{tikzpicture}
\caption{\small For $p=2$. The {\em Felli-Schneider region}, or symmetry breaking region, appears in dark grey and is defined by $a<0$ and $a<b<b_{\mathrm{FS}}(a)$. And symmetry holds in the brown region defined by $a<0$ and $ b_{\mathrm{FS}}(a)\leq b<a+1$, also $0\leq a<a_c$ and $a\leq b<a+1$.}
\label{F1}
\end{figure}
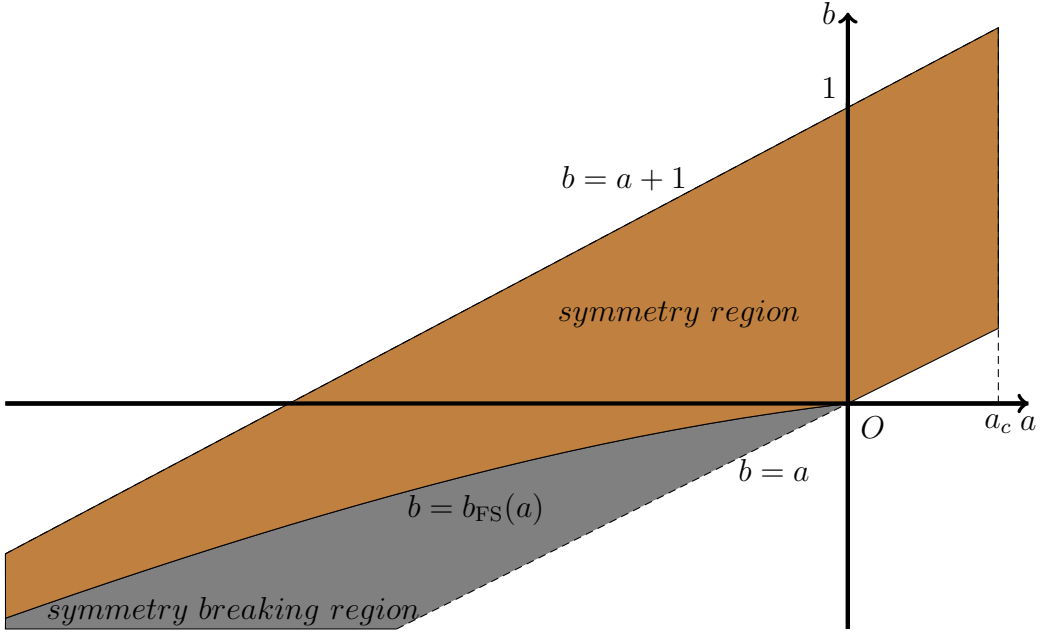

    For $p=2$, there are complete results for the above problem. When $b=a+1$ or $a<0$ and $b=a$, Catrina and Wang \cite{CW01} proved that $C_{\mathrm{CKN}}$ is not achieved and for other cases it is always achieved. Furthermore, when $a<0$ and $a< b<b_{\mathrm{FS}}(a)$, where
    \[
    b_{\mathrm{FS}}(a):=\frac{N(a_c-a)}{2\sqrt{(a_c-a)^2+N-1}}
    +a-a_c,
    \]
    Felli and Schneider \cite{FS03} proved the extremal function is non-radial. However, Dolbeault et al. \cite{DELT09} proved that when $0\leq a<a_c$ and $a\leq b<a+1$, the extremal function is symmetry. Finally, in a celebrated paper, Dolbeault, Esteban and Loss \cite{DEL16} proved an optimal rigidity result by using the so-called {\em carr\'{e} du champ} method that when $a<0$ and $b_{\mathrm{FS}}(a)\leq b<a+1$, the extremal function is symmetry, and we refer to \cite{Do21} for an overall review about this method. See the previous results shown as in Figure \ref{F1}.

    For the general case $1<p<N$, we refer to
     \cite{Ho97,SW03,TY04,WaWi00,ZZ15} for the existence and non-existence of minimizers for best constant in \eqref{cknwit}. Caldiroli and Musina \cite{CM13} obtained the symmetry breaking result, that is, if $a<b<a+1$ and
    \[
    \frac{(1+a-b)(a_c-a)}{(a_c-a+b)}
    >\sqrt{\frac{N-1}{N/(1+a-b)-1}}.
    \]
    then in \eqref{cknwit} no minimizer for best constant $C_{\mathrm{CKN}}$ is radial, see also \cite{BW02}. Furthermore,
    Lam and Lu \cite{LL17} obtained the symmetry region
    \[
    a>0\quad \mbox{and}\quad a\leq b<a+1,
    \]
    see also \cite{Ho97}.
    Recently, Ciraolo and Corso \cite{CC22} obtained partial symmetry result, that is, the authors proved the symmetry region:  if $a=b$ or $p< \frac{N}{2(1+a-b)}$,
    \[
    \frac{(1+a-b)(a_c-a)}{(a_c-a+b)}
    \leq \sqrt{\frac{N-2}{N/(1+a-b)-2}},
    \]
    furthermore they conjectured the optimal symmetry region is
    \[
    \frac{(1+a-b)(a_c-a)}{(a_c-a+b)}
    \leq \sqrt{\frac{N-1}{N/(1+a-b)-1}}.
    \]
    This is still an interesting open problem. 

    As mentioned previous, once the extremal functions of \eqref{cknwit} are well understood, it is natural to study the quantitative stability of (CKN) inequality \eqref{cknwit} by asking whether the deviation of a given function from attaining equality in \eqref{cknwit} controls its distance from the family of extremal functions.
    In the symmetry region, there are many papers concerned the stability of inequalities with potentials. R\u{a}dulescu et al. \cite{RSW02} gave the remainder terms of Hardy-Sobolev inequality for the case $p=2$ and $a=0$.
    Wang and Willem \cite{WaWi03} studied (CKN) inequalities with Lebesgue-type remainder terms, see also \cite{ACP05,DT23-1,ST18} for remainder terms of weighted Hardy inequalities.
    Wei and Wu \cite{WW22} established the stability of the profile decompositions to the (CKN) inequality \eqref{cknwit} when $p=2$ in the optimal symmetry region and also gave the gradient type remainder term for $N\geq 3$, see also our recent work \cite{DT23R2} for $N=2$. Recently, we established the stability of inequality \eqref{cknwit} when $a=0$ which extends the work of Figalli and Zhang \cite{FZ22} to Hardy case, see \cite{DT23}. Therefore, it is natural to consider the stability of inequality \eqref{cknwit} in more general symmetry origin. In present paper, we consider the case obtained by Lam and Lu \cite{LL17} and also \cite{Ho97}.

\subsection{Problem setup and main results}\label{subsectmr}

    For the consistency of the paper as in \cite{LL17}, let us perform the change of exponents
    \[
    a=\frac{\mu}{p},\quad b=\frac{s}{r},
    \]
    then \eqref{cknwit} can be rewritten in the following equivalent form:
    \begin{equation*}
    \left(\int_{\mathbb{R}^N}|u|^r\frac{\mathrm{d}x}{|x|^s}
    \right)^{\frac{p}{r}}
    \leq C\int_{\mathbb{R}^N}|\nabla u|^p\frac{\mathrm{d}x}{|x|^{\mu}},\quad \mbox{for all}\quad u\in C^\infty_c(\mathbb{R}^N).
    \end{equation*}
    Lam and Lu \cite{LL17} (see also \cite{Ho97}) proved the best constant is achieved by radial symmetry functions of above inequality in the following range:
    \begin{equation}\label{cknib1c}
    1<p<N,\quad 0<\mu<N-p,\quad \frac{\mu}{p}\leq \frac{s}{r}<\frac{\mu}{p}+1,\quad r=\frac{p(N-s)}{N-p-\mu},
    \end{equation}
    which is equivalent to $a>0$ and $a\leq b<a+1$ as in \eqref{cknwitc}.
    Define $\mathcal{D}^{1,p}_\mu(\mathbb{R}^N)$ as the closure of $C^\infty_c(\mathbb{R}^N)$
    with respect to the norm
    \[
    \|u\|_{\mathcal{D}^{1,p}_\mu(\mathbb{R}^N)}
    :=\left(\int_{\mathbb{R}^N}|x|^{-\mu}|\nabla u|^p\mathrm{d}x\right)^\frac{1}{p}.
    \]
    Define also
    \begin{equation*}
    \mathcal{S}=\mathcal{S}(N,p,\mu,s)
    :=\inf_{u\in\mathcal{D}^{1,p}_\mu(\mathbb{R}^N)\backslash\{0\}}
    \frac{\int_{\mathbb{R}^N}|\nabla u|^p\frac{\mathrm{d}x}{|x|^{\mu}}}
    {\left(\int_{\mathbb{R}^N}|u|^r\frac{\mathrm{d}x}{|x|^s}\right)
    ^{\frac{p}{r}}},
    \end{equation*}
    thus $\mathcal{S}>0$ and
    \begin{equation}\label{cknp}
    \int_{\mathbb{R}^N}|\nabla u|^p\frac{\mathrm{d}x}{|x|^{\mu}}\geq \mathcal{S}
    \left(\int_{\mathbb{R}^N}|u|^r\frac{\mathrm{d}x}{|x|^s}
    \right)^{\frac{p}{r}},
    \quad \mbox{for all}\quad u\in\mathcal{D}^{1,p}_\mu(\mathbb{R}^N).
    \end{equation}
    Note that \eqref{cknp} is equivalent to
\[
\frac{\left(\int_{\mathbb{R}^N}|\nabla u|^p\frac{\mathrm{d}x}{|x|^{\mu}}\right)^{\frac{1}{p}}}{\left(\int_{\mathbb{R}^N}|u|^r\frac{\mathrm{d}x}{|x|^s}
    \right)^{\frac{1}{r}}}
    -\mathcal{S}^{\frac{1}{p}}\geq 0,
    \quad \mbox{for all}\quad u\in \mathcal{D}^{1,p}_\mu(\mathbb{R}^N)\setminus\{0\}.
\]

    Then from \cite[Theorem 5.1]{LL17}, we directly deduce the following conclusion.

    \vskip0.25cm

    \noindent{\bf Theorem~A.} {\it Assume that \eqref{cknib1c} holds. Then the best constant $\mathcal{S}$ in \eqref{cknp} is achieved with the extremals being of the form $CU_{\lambda}$ for $C\in\mathbb{R}\setminus\{0\}$ and $\lambda>0$, where $U_{\lambda}(x)=\lambda^{\frac{N-p-\mu}{p}}U(\lambda x)$. Here,
    \begin{small}\begin{equation}\label{defeua}
    U(x)=\frac{C_{N,p,\mu,s}}
    {(1+|x|^{\frac{p-s+\mu}{p-1}})
    ^{\frac{N-p-\mu}{p-s+\mu}}}, \quad \mbox{with}\quad C_{N,p,\mu,s}=\left[(N-s)\left(\frac{N-p-\mu}{p-1}
    \right)^{p-1}\right]^{\frac{N-p-\mu}{p(p-s+\mu)}}.
    \end{equation}\end{small}
    }

\vskip0.25cm


    Note that the Euler-Lagrange equation of inequality \eqref{cknp}, up to scaling, is given by the following weighted $p$-Laplace eqution
    \begin{equation}\label{Pwha}
    -{\rm div}(|x|^{-\mu}|\nabla u|^{p-2}\nabla u)=|x|^{-s} |u|^{r-2}u\quad \mbox{in}\quad \mathbb{R}^N,\quad u\in \mathcal{D}^{1,p}_\mu(\mathbb{R}^N).
    \end{equation}
    In our recent paper \cite{DT22}, it was proved that equation \eqref{Pwha} has a unique positive radial function (up to scalings) of the form $U_\lambda$, see also \cite{Mu14,SW22}. Due to the quasi-linear property of $p$-Laplacian, let us denote $C^1_{c,0}(\mathbb{R}^N)$  be the space of compactly supported functions of class $C^1$ that are zero in a neighborhood of the origin, then we define the weighted Sobolev space $\mathcal{D}^{1,2}_{\mu,*}(\mathbb{R}^N)$ as the closure of $C^1_{c,0}(\mathbb{R}^N)$ with respect to the inner product
    \[
    \langle u,v\rangle_{\mathcal{D}^{1,2}_{\mu,*}(\mathbb{R}^N)}
    =\int_{\mathbb{R}^N}|x|^{-\mu}|\nabla U|^{p-2} \nabla u\cdot\nabla v  \mathrm{d}x,
    \]
    and the norm
    \begin{equation}\label{defd12*nb}
    \|u\|_{\mathcal{D}^{1,2}_{\mu,*}(\mathbb{R}^N)}
    :=\left(\int_{\mathbb{R}^N}|x|^{-\mu}|\nabla U|^{p-2} |\nabla u|^2  \mathrm{d}x\right)^{\frac{1}{2}}.
    \end{equation}

    \begin{remark}\label{remdefsn2}\rm
    As stated in \cite[Remark 3.1]{FZ22}, it is important for us to consider weighted that are not necessarily integrable at the origin, since $|x|^{-\mu}|\nabla U|^{p-2}\sim |x|^{\frac{(p-2)(1-s+\mu)}{p-1}-\mu}\not\in L^1(B_1)$ for $p\leq \frac{N-\mu +2(1-s+\mu)}{N+1-s}$. This is why, when defining weighted Sobolev spaces, we consider the space $C^1_{c,0}(\mathbb{R}^N)$, so that gradients vanish near zero. Of course, replacing $C^1_c(\mathbb{R}^N)$ by $C^1_{c,0}(\mathbb{R}^N)$ plays no role in the case $p> \frac{N-\mu +2(1-s+\mu)}{N+1-s}$.
    \end{remark}

    Following the work of Pistoia and Varia \cite{PV21}, our first result concerns the linearized problem related to \eqref{Pwha} at the function $U$. This leads to study the problem:
    \begin{align}\label{Ppwhla}
    & -\mathrm{div}(|x|^{-\mu}|\nabla U|^{p-2}\nabla v)-(p-2)\mathrm{div}(|x|^{-\mu}|\nabla U|^{p-4}(\nabla U\cdot\nabla v)\nabla U) \nonumber\\
    & \quad = (r-1)|x|^{-s} U^{r-2}v \quad \mbox{in}\quad \mathbb{R}^N,\quad v\in \mathcal{D}^{1,2}_{\mu,*}(\mathbb{R}^N).
    \end{align}

    It is easy to verify that $\frac{N-p-\mu}{p}U+x\cdot \nabla U$ (which equals $\frac{\partial U_{\lambda}}{\partial \lambda}|_{\lambda=1}$) solves the linear equation \eqref{Ppwhla}. We say $U$ is non-degenerate if all the solutions of \eqref{Ppwhla} result from the invariance (up to scalings) of \eqref{Pwha}. The non-degeneracy of solutions for \eqref{Pwha} is a key ingredient in analyzing the blow-up phenomena of solutions to various elliptic equations on bounded or unbounded domain in $\mathbb{R}^N$ or Riemannian manifolds whose asymptotic behavior is encoded in \eqref{defeua}, we refer to \cite{PV21} for detials.
    Therefore, it is quite natural to ask the following question:
    \begin{center}
    {\em is solution $U$ non-degenerate?}
    \end{center}
    We give an affirmative answer.

    \begin{theorem}\label{coroPpwhlpa}
    Assume that \eqref{cknib1c} holds. Then the space of solutions of \eqref{Ppwhla} has dimension one and is spanned by $(\frac{N-p-\mu}{p}U+x\cdot\nabla U)$.
    \end{theorem}

    A direct application of Theorem \ref{coroPpwhlpa} is studying the gradient stability of (CKN) inequality \eqref{cknp}. Now, let us state the main result.

    \begin{theorem}\label{thmprtp}
    Assume that \eqref{cknib1c} holds. Then there exists a constant $\mathcal{B}=\mathcal{B}(N,p,\mu,s)>0$ such that for every $u\in \mathcal{D}^{1,p}_{\mu}(\mathbb{R}^N)\setminus\{0\}$, it holds that
    \[
    \frac{\left(\int_{\mathbb{R}^N}|\nabla u|^p\frac{\mathrm{d}x}{|x|^{\mu}}\right)^{\frac{1}{p}}}{\left(\int_{\mathbb{R}^N}|u|^r\frac{\mathrm{d}x}{|x|^s}
    \right)^{\frac{1}{r}}}
    -\mathcal{S}^{\frac{1}{p}}
    \geq \mathcal{B}
    \inf_{v\in \mathcal{M}}
    \left(\frac{\|u-v\|_{\mathcal{D}^{1,p}
    _{\mu}(\mathbb{R}^N)}}{\|u\|_{\mathcal{D}^{1,p}
    _{\mu}(\mathbb{R}^N)}}\right)^\gamma,
    \]
    where $\gamma=\max\{p,2\}$ and $\mathcal{M}=\{cU_{\lambda}: c\in\mathbb{R}, \lambda>0\}$ is the set of extremal functions for (CKN) inequality \eqref{cknp}.
    \end{theorem}

    \begin{remark}\label{remp2f}\rm
    From the proof of Theorem \ref{thmprtp} as in subsection \ref{subspmr}, it is easy to verify that for $p\geq 2$,
    \begin{align*}
    \int_{\mathbb{R}^N}|\nabla u|^p\frac{\mathrm{d}x}{|x|^{\mu}}
    -\mathcal{S}
    \left(\int_{\mathbb{R}^N}|u|^r\frac{\mathrm{d}x}{|x|^s}
    \right)^{\frac{p}{r}}
    \geq C
    \inf_{v\in \mathcal{M}}
    \|u-v\|_{\mathcal{D}^{1,p}
    _{\mu}(\mathbb{R}^N)}^p, \quad \mbox{for all}\quad  u\in \mathcal{D}^{1,p}_\mu(\mathbb{R}^N),
    \end{align*}
    for some constant $C>0$, however it is false for $1<p<2$ due to Lemma \ref{lemma:rtnm2b2}.
    \end{remark}

    \begin{remark}\label{remchange}\rm
    In the proofs of our results, two different variable methods were used. In the proof of non-degenerate result, that is, for Theorem \ref{coroPpwhlpa} we will use the change $r\mapsto r^{\frac{p}{p-s+\mu}}$ in the radial case (see \cite{BCG21} for the case $p=2$) which transforms this problem into classical $p$-Laplacian linearized problem as in \cite{PV21}. However, in the proof of Theorem \ref{thmprtp}, we will follow the work of Lam and Lu \cite{LL17}, and make the change
    \begin{align*}
    v(x)=\mathcal{D}u(x):=\varrho^{-\frac{p-1}{p}}
    u(|x|^{\varrho-1}x),
    \end{align*}
    where $\varrho=\frac{N-p}{N-p-\mu}>1$, it holds that
    \begin{align*}
    \int_{\mathbb{R}^N}|u|^r\frac{\mathrm{d}x}{|x|^s}
    =\varrho\int_{\mathbb{R}^N}|v|^{r}\frac{\mathrm{d}x}{|x|^{s'}},
    \end{align*}
    where $s'=\frac{s(N-p)-N\mu}{N-p-\mu}$ (note that $r=\frac{p(N-s')}{N-p}=\frac{p(N-s)}{N-p-\mu}$ and assumption \eqref{cknib1c} implies $s'\in [0,p)$), and
    \begin{align*}
    \int_{\mathbb{R}^N}|\nabla u|^p\frac{\mathrm{d}x}{|x|^{\mu}}
    \geq \int_{\mathbb{R}^N}|\nabla v|^p \mathrm{d}x,
    \end{align*}
    furthermore, the equality holds if and only if $u$ is radial.  From above results, Lam and Lu \cite{LL17} obtained Theorem A and $\mathcal{S}=\varrho^{-(\frac{1}{r}+p-1)}\mathcal{S}'$,  where $\mathcal{S}'$ is sharp constant of Hardy-Sobolev inequality (see \cite{GY00}).
    These results will be used helpful in the proof of Lemma \ref{propcetl} which is a crucial compact theory for proving the case $1<p\leq \frac{2(N-\mu)}{N+2-s}$.
    It is worth mentioning our recent work \cite{DT23} that there is constant $c>0$ such that for all $v\in \mathcal{D}^{1,p}_{0}(\mathbb{R}^N)\setminus\{0\}$,
    \[
    \frac{\left(\int_{\mathbb{R}^N}|\nabla v|^p\mathrm{d}x\right)^{\frac{1}{p}}}
    {\left(\int_{\mathbb{R}^N}|v|^{r}\frac{\mathrm{d}x}{|x|^{s'}}\right)^{\frac{1}{r}}}
    -(\mathcal{S'})^{\frac{1}{p}}
    \geq
    c\inf_{W\in \mathcal{M}'}
    \left(\frac{\|v-W\|_{\mathcal{D}^{1,p}
    _0(\mathbb{R}^N)}}{\|v\|_{\mathcal{D}^{1,p}
    _{0}(\mathbb{R}^N)}}\right)^{\max\{2,p\}},
    \]
    where $\mathcal{M}'$ is the set of extremal functions for Hardy-Sobolev inequality.
    But it can not deduce Theorem \ref{thmprtp} directly, since under previous change of variable (also $W(x)=\mathcal{D}U(x)$),
    \[
    \|u-U\|
    _{\mathcal{D}^{1,p}_{\mu}(\mathbb{R}^N)}
    \geq
    \|v-W\|
    _{\mathcal{D}^{1,p}_{0}(\mathbb{R}^N)}.
    \]
    We will follow the arguments as those in \cite{FZ22} and also \cite{DT23} to prove our result.
    \end{remark}

\subsection{Structure of the paper}\label{subsect:structrue}

In Section \ref{sectce}, we prove some compactness results. Section \ref{sectndr} is devoted to proving the non-degeneracy of extremal function $U$. Then in Section \ref{sectspana}, we give the spectral analysis with the help of compactness. In Section \ref{sectpromr}, we first prove a local variant of Theorem \ref{thmprtp} then complete the proof of Theorem \ref{thmprtp} by Lions' concentration and compactness principle.  Finally, we collect two crucial technical estimates in Appendix \ref{sectpls}.

\noindent{\bfseries Notations.}
    Throughout this paper, we denote $B_\rho:=B(\mathbf{0},\rho)$ be the ball with radius $\rho$ centered at the origin.
Moreover, $c$, $C$, $C'$ and $C_i$ are indiscriminately used to denote various absolutely positive constants from line to line. $a\sim b$ means that $C'b\leq a\leq Cb$.

\section{\bfseries Compact embedding theorem}\label{sectce}

    Firstly, let us establish the following Poincar\'{e} type inequalities which will be useful later.
    \begin{lemma}\label{lemcztj}
    Assume that \eqref{cknib1c} holds. Then there exists $C=C(N,p,\mu,s)>0$ such that for any $\varphi\in \mathcal{D}^{1,2}_{\mu,*}(\mathbb{R}^N)$, it holds that
    \begin{align}\label{continouseb}
    \int_{\mathbb{R}^N}|x|^{-\mu}|\nabla U|^{p-2}|\nabla \varphi|^2\mathrm{d}x
    \geq C \int_{\mathbb{R}^N}
    |x|^{-s}|U|^{r-2}\varphi^2\mathrm{d}x.
    \end{align}
    Also, there exists $\vartheta=\vartheta(N,p,\mu,s)>0$ such that, for any $\rho\in (0,1)$, we have
    \begin{align}\label{continouseb1}
    \int_{\mathbb{R}^N}|x|^{-\mu}|\nabla U|^{p-2}|\nabla \varphi|^2\mathrm{d}x
    \geq \frac{C}{\rho^{\vartheta}}\int_{B_\rho}
    |x|^{-s}|U|^{r-2}\varphi^2\mathrm{d}x,
    \end{align}
    and
    \begin{align}\label{continouseb2}
    \int_{\mathbb{R}^N}|x|^{-\mu}|\nabla U|^{p-2}|\nabla \varphi|^2\mathrm{d}x
    \geq C|\log \rho|^2\int_{\mathbb{R}^N\setminus B_{\frac{1}{\rho}}}
    |x|^{-s}|U|^{r-2}\varphi^2\mathrm{d}x.
    \end{align}
    \end{lemma}

    \begin{proof}
    To prove \eqref{continouseb}, we can assume by approximation that $\varphi\in C^1_{c,0}(\mathbb{R}^N)$ as in \cite[Lemma 4.2]{FN19}.  Define
    \[
    \mathfrak{F}(u):=\int_{\mathbb{R}^N}|x|^{-\mu}|\nabla u|^p \mathrm{d}x
    -\mathcal{S}\left(\int_{\mathbb{R}^N}
    |x|^{-s}|u|^{r} \mathrm{d}x\right)^\frac{p}{r}.
    \]
    We know $U$ is a local minimum of the functional $\mathfrak{F}$, then
    \begin{align*}
    0\leq & \frac{\mathrm{d}^2}{\mathrm{d}\varepsilon^2}
    \Big|_{\varepsilon=0} \mathfrak{F}(U+\varepsilon\varphi)
    \\ = & p \int_{\mathbb{R}^N}|x|^{-\mu}|\nabla U|^{p-2}|\nabla \varphi|^2\mathrm{d}x
    + p(p-2) \int_{\mathbb{R}^N}|x|^{-\mu}|\nabla U|^{p-4}(\nabla U\cdot\nabla\varphi)^2\mathrm{d}x
    \\ & - \mathcal{S}
    \Bigg[
    p\left(p-r\right)
    \left(\int_{\mathbb{R}^N}|x|^{-s}|U|^{r} \mathrm{d}x\right)^{\frac{p}{r}-2}
    \left(\int_{\mathbb{R}^N}
    |x|^{-s}|U|^{r-2}U\varphi\mathrm{d}x\right)^2
    \\ & \quad\quad + p(r-1)
    \left(\int_{\mathbb{R}^N}|x|^{-s}|U|^{r} \mathrm{d}x\right)^{\frac{p}{r}-1}
    \int_{\mathbb{R}^N}
    |x|^{-s}|U|^{r-2}\varphi^2\mathrm{d}x
    \Bigg].
    \end{align*}
    Noting that
    \begin{align*}
    \left(\int_{\mathbb{R}^N}
    |x|^{-s}|U|^{r-2}U\varphi\mathrm{d}x\right)^2
    \leq \int_{\mathbb{R}^N}
    |x|^{-s}|U|^{r}\mathrm{d}x
    \int_{\mathbb{R}^N}
    |x|^{-s}|U|^{r-2}\varphi^2\mathrm{d}x,
    \end{align*}
    which implies that 
    \begin{align*}
    \int_{\mathbb{R}^N}|x|^{-s}|U|^{r-2}\varphi^2\mathrm{d}x
    \leq & \int_{\mathbb{R}^N}|x|^{-\mu}|\nabla U|^{p-2}|\nabla \varphi|^2\mathrm{d}x,
    \end{align*}
    if $p\geq 2$, and 
    \begin{align*}
    (p-1)\int_{\mathbb{R}^N}|x|^{-s}|U|^{r-2}\varphi^2\mathrm{d}x
    \leq & \int_{\mathbb{R}^N}|x|^{-\mu}|\nabla U|^{p-2}|\nabla \varphi|^2\mathrm{d}x,
    \end{align*}
    if $1<p<2$, due to
    \[
    \int_{\mathbb{R}^N}|x|^{-s}|U|^{r} \mathrm{d}x
    =\mathcal{S}^{\frac{r}{r-p}}.
    \]
    Thus \eqref{continouseb} holds.

    To prove \eqref{continouseb1}, we can also assume by approximation that $\varphi\in C^1_{c,0}(\mathbb{R}^N)$ and we apply the Sobolev inequality with radial weights. More precisely, for any $\rho\in (0,1)$ we have
    \begin{align}\label{wsit}
    \int_{\mathbb{R}^N}|x|^{-\mu}|\nabla U|^{p-2}|\nabla \varphi|^2\mathrm{d}x
    \geq & c\left(\int_{B_\rho}
    |x|^{-\frac{s\zeta}{2}}|\varphi|^\zeta\mathrm{d}x
    \right)^{\frac{2}{\zeta}},
    \end{align}
    for some $\zeta>2$.
    In fact, since $|x|^{-\mu}|\nabla U|^{p-2}\sim |x|^{\frac{(p-2)(1-s+\mu)}{p-1}-\mu}$ inside $B_\rho$, then applying \eqref{continouseb} and the Sobolev inequality with radial weights (see \cite[Section 2.1]{Ma85}) we deduce
    \begin{align*}
    \int_{\mathbb{R}^N}|x|^{-\mu}|\nabla U|^{p-2}|\nabla \varphi|^2\mathrm{d}x
    \geq & c_1\int_{\mathbb{R}^N}(|x|^{-s}|U|^{r-2}\varphi^2
    +|x|^{-\mu}|\nabla U|^{p-2}|\nabla \varphi|^2)\mathrm{d}x
    \\ \geq & c_2\int_{B_\rho}
    |x|^{-s}\left(\varphi^2
    +|x|^{\frac{(p-2)(1-s+\mu)}{p-1}-\mu+s}|\nabla \varphi|^2\right)\mathrm{d}x
    \\ \geq & c_2\int_{B_\rho}
    |x|^{-s}\left(\varphi^2+|x|^{2}|\nabla \varphi|^2\right)\mathrm{d}x
    \\ \geq & c_3\int_{B_\rho}
    \left(\left(|x|^{-\frac{s}{2}}\varphi\right)^2+|x|^{2}\left|\nabla \left(|x|^{-\frac{s}{2}}\varphi\right)\right|^2\right)\mathrm{d}x
    \\ \geq & c_4\left(\int_{B_\rho}
    |x|^{-\frac{s\zeta}{2}}|\varphi|^\zeta\mathrm{d}x
    \right)^{\frac{2}{\zeta}},
    \end{align*}
    due to $\frac{(p-2)(1-s+\mu)}{p-1}-\mu+s \leq 2$ (note that assumption \eqref{cknib1c} implies $s-\mu<p$),
    where $\zeta=\zeta(N)>2$.
    Therefore, by H\"{o}lder inequality we obtain
    \begin{align*}
    \int_{B_\rho}
    |x|^{-s}|U|^{p-2}\varphi^2\mathrm{d}x
    \leq & C_1 \int_{B_\rho}
    |x|^{-s}\varphi^2\mathrm{d}x
    \\ \leq & C_2 \rho^{N(1-\frac{\zeta}{2})}\left(\int_{B_\rho}
    |x|^{-\frac{s\zeta}{2}}|\varphi|^\zeta\mathrm{d}x
    \right)^{\frac{2}{\zeta}}
    \\ \leq & C_3 \rho^{N(1-\frac{\zeta}{2})}\int_{\mathbb{R}^N}
    |x|^{-\mu}|\nabla U|^{p-2}|\nabla \varphi|^2\mathrm{d}x,
    \end{align*}
    as desired for $\vartheta=N(\frac{\zeta}{2}-1)>0$ in \eqref{continouseb1}.

    To prove \eqref{continouseb2}, for $\rho\in (0,1)$ we define
    \begin{eqnarray*}
    \eta_\rho(x):=
    \left\{ \arraycolsep=1.5pt
       \begin{array}{ll}
        0,\ \ &{\rm for}\ \ |x|<\rho^{-\frac{1}{2}},\\[2mm]
        \frac{2\log |x|-|\log \rho|}{|\log \rho|},\ \ &{\rm for}\ \  \rho^{-\frac{1}{2}}\leq |x|\leq \rho^{-1},\\[2mm]
        1,\ \ &{\rm for}\ \ |x|> \rho^{-1},
        \end{array}
    \right.
    \end{eqnarray*}
    and $\phi_\rho:=\eta_\rho\varphi$.
    We note that, since
    \begin{equation*}
    U(\tau)\sim
    (1+\tau^{\frac{p-s+\mu}{p-1}})
    ^{-\frac{N-p-\mu}{p-s+\mu}},\quad
    U'(\tau)\sim
    (1+\tau^{\frac{p-s+\mu}{p-1}})
    ^{-\frac{N-s}{p-s+\mu}}\tau^{\frac{1-s+\mu}{p-1}},
    \end{equation*}
    then
    \begin{equation*}
    |U|^{r-2}\sim
    (1+\tau^{\frac{p-\beta}{p-1}})
    ^{-\frac{(N-s)(p-2)}{p-s+\mu}-2},\quad
    |\nabla U|^{p-2}\sim
    (1+\tau^{\frac{p-s+\mu}{p-1}})
    ^{-\frac{(N-s)(p-2)}{p-s+\mu}}\tau^{\frac{(1-s+\mu)(p-2)}{p-1}}.
    \end{equation*}
    Thanks to Fubini's theorem and using polar coordinates,
    \begin{small}\begin{align*}
    \int_{\mathbb{R}^N}
    |x|^{-s}|U|^{r-2}|\phi_\rho|^2\mathrm{d}x
    \leq & C_1\int_{\mathbb{S}^{N-1}}\int^{+\infty}_0    \tau^{N-1-s}(1+\tau^{\frac{p-s+\mu}{p-1}})
    ^{-\frac{(N-s)(p-2)}{p-s+\mu}-2}
    |\phi_\rho(\tau\theta)|^2
    \mathrm{d}\tau\mathrm{d}\theta
    \\ \leq & C_2\int_{\mathbb{S}^{N-1}}\int^{+\infty}_0    \tau^{N-1-s}(1+\tau^{\frac{p-s+\mu}{p-1}})
    ^{-\frac{(N-s)(p-2)}{p-s+\mu}-2}
    \\ & \quad \times
    \int^{+\infty}_\tau |\phi_\rho(\sigma\theta)||\nabla\phi_\rho(\sigma\theta)|
    \mathrm{d}\sigma\mathrm{d}\tau\mathrm{d}\theta
    \\ \leq & C_3\int_{\mathbb{S}^{N-1}}\int^{+\infty}_0 \int^{\sigma}_0    \tau^{N-1-s}(1+\tau^{\frac{p-s+\mu}{p-1}})
    ^{-\frac{(N-s)(p-2)}{p-s+\mu}-2}
    \\ & \quad \times
    |\phi_\rho(\sigma\theta)||\nabla\phi_\rho(\sigma\theta)|
    \mathrm{d}\tau\mathrm{d}\sigma\mathrm{d}\theta
    \\ \leq & C_4\int_{\mathbb{S}^{N-1}}\int^{+\infty}_0     \sigma^{N-s}(1+\sigma^{\frac{p-s+\mu}{p-1}})
    ^{-\frac{(N-s)(p-2)}{p-s+\mu}-2}
    |\phi_\rho(\sigma\theta)||\nabla\phi_\rho(\sigma\theta)|
    \mathrm{d}\sigma\mathrm{d}\theta.
    \end{align*}\end{small}
    Thus, by Cauchy-Schwarz inequality we obtain
    \begin{small}\begin{align*}
    \int_{\mathbb{R}^N}
    |x|^{-s}|U|^{r-2}|\phi_\rho|^2\mathrm{d}x
    \leq & C_5\left(\int_{\mathbb{S}^{N-1}}\int^{+\infty}_0     \sigma^{N+1-s}(1+\sigma^{\frac{p-s+\mu}{p-1}})
    ^{-\frac{(N-s)(p-2)}{p-s+\mu}-2}
    |\nabla\phi_\rho(\sigma\theta)|
    \mathrm{d}\theta\right)^{\frac{1}{2}}
    \\ & \quad \times
    \left(\int_{\mathbb{S}^{N-1}}\int^{+\infty}_0     \sigma^{N-1-s}(1+\sigma^{\frac{p-s+\mu}{p-1}})
    ^{-\frac{(N-s)(p-2)}{p-s+\mu}-2}
    |\phi_\rho(\sigma\theta)|
    \mathrm{d}\sigma\right)^{\frac{1}{2}},
    \end{align*}\end{small}
    and since the last term in the right hand coincides with $\|\phi_\rho\|_{L^2_{s,*}(\mathbb{R}^N)}$ (up to a multiplicative constant), we conclude that
    \begin{align}\label{pfleql}
    \int_{\mathbb{R}^N}
    |x|^{-s}|U|^{r-2}|\phi_\rho|^2\mathrm{d}x
    \leq & C_6\int_{\mathbb{S}^{N-1}}\int^{+\infty}_0     \sigma^{N+1-s}(1+\sigma^{\frac{p-s+\mu}{p-1}})
    ^{-\frac{(N-s)(p-2)}{p-s+\mu}-2}
    |\nabla\phi_\rho(\sigma\theta)|
    \mathrm{d}\theta
    \nonumber\\ \leq & C_7 \int_{\mathbb{R}^N}|x|^{2-s}(1+|x|^{\frac{p-s+\mu}{p-1}})
    ^{-\frac{(N-s)(p-2)}{p-s+\mu}-2}|\nabla\phi_\rho(x)|^2
    \mathrm{d}x.
    \end{align}
    Then we have
    \begin{align*}
    \int_{\mathbb{R}^N\setminus B_{\frac{1}{\rho}}}
    |x|^{-s}|U|^{r-2}|\phi_\rho|^2\mathrm{d}x
    \leq & C_1\int_{\mathbb{R}^N}
    |x|^{-s}|U|^{r-2}|\phi_\rho|^2\mathrm{d}x
    \\ \leq & C_2 \int_{\mathbb{R}^N}|x|^{2-s}(1+|x|^{\frac{p-s+\mu}{p-1}})
    ^{-\frac{(N-s)(p-2)}{p-s+\mu}-2}|\nabla\phi_\rho(x)|^2
    \mathrm{d}x
    \\ \leq & C_3 \int_{\mathbb{R}^N}|x|^{2-s}(1+|x|^{\frac{p-s+\mu}{p-1}})
    ^{-\frac{(N-s)(p-2)}{p-s+\mu}-2}\eta_\rho^2|\nabla\varphi|^2
    \mathrm{d}x
    \\ & + C_3 \int_{\mathbb{R}^N}|x|^{2-s}(1+|x|^{\frac{p-s+\mu}{p-1}})
    ^{-\frac{(N-s)(p-2)}{p-s+\mu}-2}|\nabla\eta_\rho|^2
    \varphi^2\mathrm{d}x
    \\ \leq & C_3 \int_{\mathbb{R}^N\setminus B_{\frac{1}{\rho^{1/2}}}}|x|^{2-s}(1+|x|^{\frac{p-s+\mu}{p-1}})
    ^{-\frac{(N-s)(p-2)}{p-s+\mu}-2}|\nabla\varphi|^2
    \mathrm{d}x
    \\ & + 4C_3 |\log \rho|^{-2}
    \int_{B_{\frac{1}{\rho}}\setminus B_{\frac{1}{\rho^{1/2}}}}|x|^{-s}(1+|x|^{\frac{p-s+\mu}{p-1}})
    ^{-\frac{(N-s)(p-2)}{p-s+\mu}-2}\varphi^2
    \mathrm{d}x
    \\ \leq & C_4 \rho^{\frac{p-s+\mu}{2(p-1)}}\int_{\mathbb{R}^N\setminus B_{\frac{1}{\rho^{1/2}}}}|x|^{-\mu}|\nabla U|^{p-2}|\nabla\varphi|^2
    \mathrm{d}x
    \\ & + C_4 |\log \rho|^{-2}\int_{B_{\frac{1}{\rho}}\setminus B_{\frac{1}{\rho^{1/2}}}}|x|^{-s}
    |U|^{r-2}\varphi^2\mathrm{d}x.
    \end{align*}
    Thus, combining with \eqref{continouseb} we deduce that
    \begin{align*}
    \int_{\mathbb{R}^N\setminus B_{\frac{1}{\rho}}}
    |x|^{-s}|U|^{r-2}\varphi^2\mathrm{d}x
    \leq &   C_4 \rho^{\frac{p-s+\mu}{2(p-1)}}\int_{\mathbb{R}^N}
    |x|^{-\mu}|\nabla U|^{p-2}|\nabla\varphi|^2
    \mathrm{d}x
    \\ & + C_4 |\log \rho|^{-2}\int_{\mathbb{R}^N}|x|^{-s}
    |U|^{r-2}\varphi^2\mathrm{d}x
    \\ \leq & C_5 |\log \rho|^{-2}\int_{\mathbb{R}^N}|x|^{-\mu}|\nabla U|^{p-2}|\nabla\varphi|^2
    \mathrm{d}x,
    \end{align*}
    due to $\frac{p-s+\mu}{2(p-1)}>0$ and $\rho\in (0,1)$, \eqref{continouseb2} holds.
    \end{proof}

    Now, let us give a crucial compact embedding theorem as follows.
    \begin{proposition}\label{propcet}
    Assume that \eqref{cknib1c} holds. The space $\mathcal{D}^{1,2}_{\mu,*}(\mathbb{R}^N)$ compactly embeds into $L^2_{s,*}(\mathbb{R}^N)$, where weighted Sobolev space $L^2_{s,*}(\mathbb{R}^N)$ is defined as the set of measurable functions with respect to the norm
    \begin{equation}\label{defd12*na}
    \|\varphi\|_{L^{2}_{s,*}(\mathbb{R}^N)}:
    =\left(\int_{\mathbb{R}^N}
    |x|^{-s} U^{r-2}\varphi^2\mathrm{d}x\right)^{\frac{1}{2}}.
    \end{equation}
    \end{proposition}

    \begin{proof}
    Let $\{\varphi_n\}$ be a sequence of functions in $\mathcal{D}^{1,2}_{\mu,*}(\mathbb{R}^N)$ with uniformly bounded norm. It follows from \eqref{continouseb} that $\|\varphi_n\|_{L^2_{s,*}(\mathbb{R}^N)}$ are uniformly bounded as well.

    Since both $|x|^{-\mu}|\nabla U|^{p-2}$ and $|x|^{-s}|U|^{r-2}$ are locally bounded away from zero and infinity in $\mathbb{R}^N\setminus\{\mathbf{0}\}$, by Rellich-Kondrachov Theorem and a diagonal argument we deduce that, up to a subsequence, there exists $\varphi\in \mathcal{D}^{1,2}_{\mu,*}(\mathbb{R}^N)\cap L^2_{s,*}(\mathbb{R}^N)$ such that $\varphi_n\rightharpoonup \varphi$ in $\mathcal{D}^{1,2}_{\mu,*}(\mathbb{R}^N)\cap L^2_{s,*}(\mathbb{R}^N)$ and $\varphi_n\to \varphi$ locally in $L^2_{s,*}(\mathbb{R}^N\setminus\{\mathbf{0}\})$.

    Also, it follows from \eqref{continouseb1} and \eqref{continouseb2} that, there exists $\vartheta>0$ such that for any $\rho\in (0,1)$,
    \begin{align*}
    \int_{B_\rho}
    |x|^{-s}|U|^{r-2}\varphi^2\mathrm{d}x
    \leq C\rho^{\vartheta},
    \quad
    \int_{B_{\frac{1}{\rho}}}
    |x|^{-s}|U|^{r-2}\varphi^2\mathrm{d}x
    \leq \frac{C}{|\log \rho|^2}.
    \end{align*}
    We conclude the proof by defining the compact set $K_\rho:=\overline{B_{\frac{1}{\rho}}}\setminus B_\rho$ and applying the strong convergence of $\varphi_n$ in $K_\rho$, together with the arbitrariness of $\rho$ (that can be chosen arbitrarily small).
    \end{proof}

    As mentioned in \cite{FZ22}, because of the crucial inequality \ref{uinx2pl}, we shall see that Proposition \ref{propcet} allows us to deal with the case $\frac{2(N-\mu)}{N+2-s}< p<2$ when we show the stability of (CKN) inequality \eqref{cknp}. However, when $1<p\leq\frac{2(N-\mu)}{N+2-s}$ which implies $p<r\leq 2$, we will need a much more delicate compactness result that we now present.

    \begin{lemma}\label{propcetl}
    Assume that \eqref{cknib1c} holds and $1<p\leq \frac{2(N-\mu)}{N+2-s}$. Let $\{v_n\}$ be a sequence of functions in $\mathcal{D}^{1,p}_{\mu}(\mathbb{R}^N)$ satisfying
    \begin{align}\label{propcetll}
    & \int_{\mathbb{R}^N}
    |x|^{-\mu}(|\nabla U|+\varepsilon_n|\nabla v_n|)^{p-2}|\nabla v_n|^2\mathrm{d}x\leq 1,
    \end{align}
    where $\varepsilon_n\in (0,1)$ is a sequence of positive numbers converging to $0$. Then, up to a subsequence, $v_n$ convergence weakly in $\mathcal{D}^{1,p}_{\mu}(\mathbb{R}^N)$ to some $v\in\mathcal{D}^{1,p}_{\mu}(\mathbb{R}^N)\cap L^2_{s,*}(\mathbb{R}^N)$. Also, given any constant $\mathcal{C}\geq 0$ it holds
    \begin{align}\label{propcetlc}
    \int_{\mathbb{R}^N}|x|^{-s}
    \frac{(U+\mathcal{C}|\varepsilon_nv_n|)^{r}}
    {U^2+|\varepsilon_nv_n|^2}|v_n|^2\mathrm{d}x
    \to \int_{\mathbb{R}^N}|x|^{-s}
    U^{r-2}|v|^2\mathrm{d}x.
    \end{align}
    \end{lemma}

    \begin{proof}
    Observe that, since $1<p<2$, by H\"{o}lder inequality we have
    \begin{small}
    \begin{align*}
    \int_{\mathbb{R}^N}|x|^{-\mu}|\nabla v_n|^p\mathrm{d}x
    \leq & \left(\int_{\mathbb{R}^N}|x|^{-\mu}
    (|\nabla U|+\varepsilon_n|\nabla v_n|)^{p-2}|\nabla v_n|^2\mathrm{d}x \right)^{\frac{p}{2}}
    \\ &\times
    \left(\int_{\mathbb{R}^N}|x|^{-\mu}
    (|\nabla U|+\varepsilon_n|\nabla v_n|)^{p}dx \right)^{1-\frac{p}{2}} \\
    \leq & C(N,p,\mu,s)\left(\int_{\mathbb{R}^N}|x|^{-\mu}
    (|\nabla U|+\varepsilon_n|\nabla v_n|)^{p-2}|\nabla v_n|^2\mathrm{d}x \right)^{\frac{p}{2}}
    \\ &\times
    \left(1+\varepsilon_n^p\int_{\mathbb{R}^N}|x|^{-\mu}|\nabla v_n|^{p}\mathrm{d}x \right)^{1-\frac{p}{2}},
    \end{align*}
    \end{small}
    that combined with \eqref{propcetll} gives
    \begin{align*}
    \left(\int_{\mathbb{R}^N}|x|^{-\mu}|\nabla v_n|^p\mathrm{d}x\right)^{\frac{2}{p}}
    \leq & C\int_{\mathbb{R}^N}|x|^{-\mu}
    (|\nabla U|+\varepsilon_n|\nabla v_n|)^{p-2}|\nabla v_n|^2\mathrm{d}x
    \leq C.
    \end{align*}
    Furthermore, \eqref{continouseb} indicates
    \begin{align*}
    \int_{\mathbb{R}^N}
    |x|^{-s}|U|^{r-2}v_n^2\mathrm{d}x
    \leq & C\int_{\mathbb{R}^N}|x|^{-\mu}|\nabla U|^{p-2}|\nabla v_n|^2\mathrm{d}x
    \\
    \leq & C\int_{\mathbb{R}^N}|x|^{-\mu}
    (|\nabla U|+\varepsilon_n|\nabla v_n|)^{p-2}|\nabla v_n|^2\mathrm{d}x
    \leq C.
    \end{align*}
    Thus, up to a subsequence, $v_n$ convergence weakly in $\mathcal{D}^{1,p}_{\mu}(\mathbb{R}^N)$ to some $v\in\mathcal{D}^{1,p}_{\mu}(\mathbb{R}^N)\cap L^2_{s,*}(\mathbb{R}^N)$. Hence, to conclude the proof, we need to show the validity of \eqref{propcetlc}.

    Following the work of Lam and Lu \cite{LL17} shown as in Remark \ref{remchange}, we make the change
    \begin{align}\label{changey}
    w_n(x)=\mathcal{D}v_n(x):=\left(\frac{1}{\varrho}\right)^{\frac{p-1}{p}}
    v_n(|x|^{\varrho-1}x),
    \end{align}
    where $\varrho=\frac{N-p}{N-p-\mu}>1$, it holds that
    \begin{align}\label{changel}
    |\nabla w_n(x)|\leq \varrho^{\frac{1}{p}}|x|^{\varrho-1}|\nabla v_n(|x|^{\varrho-1}x)|,
    \end{align}
    and
    \begin{align*}
    \int_{\mathbb{R}^N}|x|^{-s}
    \frac{(U+\mathcal{C}|\varepsilon_nv_n|)^{r}}
    {U^2+|\varepsilon_nv_n|^2}|v_n|^2\mathrm{d}x
    = \varrho \int_{\mathbb{R}^N}|x|^{-[N+(s-N)\varrho]}
    \frac{(W+\mathcal{C}|\varepsilon_nw_n|)^{r}}
    {W^2+|\varepsilon_nw_n|^2}|w_n|^2\mathrm{d}x,
    \end{align*}
    where $W(x)=\mathcal{D}U(x)$ and $0\leq N+(s-N)\varrho<p$, see \cite[Lemma 2.2]{LL17}. Also,
    \begin{align*}
    \int_{\mathbb{R}^N}|x|^{-s}
    U^{r-2}|v|^2\mathrm{d}x
    =\varrho \int_{\mathbb{R}^N}|x|^{-[N+(s-N)\varrho]}
    U^{r-2}|w|^2\mathrm{d}x,
    \end{align*}
    where $w(x)=\mathcal{D}v(x)$. From our recent work \cite[Lemma 2.2]{DT23}, it holds that
    \begin{align}\label{changccj}
    \int_{\mathbb{R}^N}|x|^{-[N+(s-N)\varrho]}
    \frac{(W+\mathcal{C}|\varepsilon_nw_n|)^{r}}
    {W^2+|\varepsilon_nw_n|^2}|w_n|^2\mathrm{d}x
    \to \int_{\mathbb{R}^N}|x|^{-[N+(s-N)\varrho]}
    W^{r-2}|w|^2\mathrm{d}x,
    \end{align}
    for $0<N+(s-N)\varrho<p$, if
    \begin{align}\label{changcc}
    \int_{\mathbb{R}^N}
    (|\nabla W|+\varepsilon_n|\nabla w_n|)^{p-2}|\nabla w_n|^2\mathrm{d}x\leq 1.
    \end{align}
    Otherwise, if $N+(s-N)\varrho$ which implies $r=p^*$, then \cite[Lemma 3.4]{FZ22} also indicates \eqref{changccj} if \eqref{changcc} holds.
    So we only need to prove \eqref{changcc} under the change \eqref{changey}. Note that the function $f(t)=(a+bt)^{p-2}t^2$ where $a\geq 0$ and $0<b<1$, is nondecreasing for $t\geq 0$, then from \eqref{changel} we obtain
    \begin{align*}
    & \int_{\mathbb{R}^N}
    |x|^{-\mu}(|\nabla U|+\varepsilon_n|\nabla v_n|)^{p-2}|\nabla v_n|^2\mathrm{d}x
    \\
    & \quad = \int_{\mathbb{R}^N}
    \frac{\left(|\nabla U(|x|^{\varrho-1}x)|+\varepsilon_n|\nabla v_n(|x|^{\varrho-1}x)|\right)^{p-2}|\nabla v_n(|x|^{\varrho-1}x)|^2}
    {||x|^{\varrho-1}x|^{\mu}}
    \varrho|x|^{N(\varrho-1)}\mathrm{d}x
    \\
    & \quad \geq \frac{1}{\varrho}\int_{\mathbb{R}^N}
    \frac{(|\nabla W|+\varepsilon_n|\nabla w_n|)^{p-2}|\nabla w_n|^2}
    {|x|^{p(\varrho-1)}||x|^{\varrho-1}x|^{\mu}}
    \varrho|x|^{N(\varrho-1)}\mathrm{d}x
    \\
    & \quad = \int_{\mathbb{R}^N}
    (|\nabla W|+\varepsilon_n|\nabla w_n|)^{p-2}|\nabla w_n|^2\mathrm{d}x.
    \end{align*}
    Therefore, from \eqref{propcetll} we deduce \eqref{changcc} holds. The proof is completed.
    \end{proof}

    An important consequence of Lemma \ref{propcetl} is the following weighted Orlicz-type Poincar\'{e} inequality:

    \begin{corollary}\label{propcetlpi}
    Assume that \eqref{cknib1c} holds and $1<p\leq \frac{2(N-\mu)}{N+2-s}$. There exists $\varepsilon_0>0$ small such that the following holds:
    For any $\varepsilon \in (0,\varepsilon_0)$ and any radial function $v\in \mathcal{D}^{1,p}_{\mu}(\mathbb{R}^N)\cap \mathcal{D}^{1,2}_{\mu,*}(\mathbb{R}^N)$ satisfying
    \begin{align*}
    & \int_{\mathbb{R}^N}|x|^{-\mu}(|\nabla U|+\varepsilon|\nabla v|)^{p-2}|\nabla v|^2\mathrm{d}x\leq 1,
    \end{align*}
    we have
    \begin{align}\label{propcetlcpi}
    \int_{\mathbb{R}^N}|x|^{-s}
    (U+ \varepsilon |v| )^{r-2}|v|^2\mathrm{d}x
    \leq C\int_{\mathbb{R}^N}|x|^{-\mu}(|\nabla U|+\varepsilon|\nabla v|)^{p-2}|\nabla v|^2\mathrm{d}x.
    \end{align}
    \end{corollary}

    \begin{proof}
    Making the change
    \[
    w(x)=\mathcal{D}v(x):=\left(\frac{1}{\varrho}\right)^{\frac{p-1}{p}}
    v(|x|^{\varrho-1}x),
    \]
    and also $W(x)=\mathcal{D}U(x)$,
    then based on Lemma \ref{propcetl}, this proof can be deduced directly from the proof of
    \cite[Corollary 3.5]{FZ22} with minor changes, so we omit it.
    \end{proof}

\section{{\bfseries Non-degenerate result}}\label{sectndr}

    First of all, a straightforward computation shows that (\ref{Ppwhla}) is equivalent to
    \begin{small}
    \begin{align}\label{rwlp}
    & -|x|^2\Delta v +\mu (x\cdot\nabla v)
    -(p-2) \sum^{N}_{i,j=1}\frac{\partial^2 v}{\partial x_i\partial x_j}x_i x_j
    -\frac{(p-2)(N-s)}{1+|x|^{\frac{p-s+\mu}{p-1}}} (x\cdot\nabla v) \nonumber\\
    & \quad = (r-1)C_{N,p,\mu,s}^{r-p}\left(\frac{N-p-\mu}{p-1}\right)^{2-p}
    \frac{|x|^{\frac{p-s+\mu}{p-1}}}{(1+|x|^{\frac{p-s+\mu}{p-1}})^2}v
    \quad \mbox{in}\quad \mathbb{R}^N,\quad v\in \mathcal{D}^{1,2}_{\mu,*}(\mathbb{R}^N).
    \end{align}
    \end{small}
    
    Then by using the standard spherical decomposition and making the change of variable $r\mapsto r^{\frac{p}{p-s+\mu}}$, we can characterize all solutions to the linearized problem (\ref{rwlp}). That is, we have the non-degenerate conclusion shown as in Theorem \ref{coroPpwhlpa}.

    \subsection{Proof of Theorem \ref{coroPpwhlpa}}

    Since $U$ is radial we can make a partial wave decomposition of (\ref{rwlp}), namely
    \begin{equation}\label{Ppwhl2defvdp}
    v(x)=v(\rho,\theta)=\sum^{\infty}_{k=0}\sum^{l_k}_{m=1}
    \varphi_{k,m}(\rho)\Psi_{k,m}(\theta),
    \end{equation}
    where $\rho=|x|$, $\theta=\frac{x}{|x|}\in \mathbb{S}^{N-1}$, and
    \begin{equation*}
    \varphi_{k,m}(\rho)=
    \int_{\mathbb{S}^{N-1}}v(\rho,\theta)\Psi_{k,m}(\theta)d\theta.
    \end{equation*}
    Here $\Psi_{k,m}(\theta)$ denotes the $k$-th spherical harmonic, i.e., it satisfies
    \begin{equation}\label{deflk}
    -\Delta_{\mathbb{S}^{N-1}}\Psi_{k,m}=\lambda_k \Psi_{k,m},
    \end{equation}
    where $\Delta_{\mathbb{S}^{N-1}}$ is the Laplace-Beltrami operator on $\mathbb{S}^{N-1}$ with the standard metric and  $\lambda_k$ is the $k$-th eigenvalue of $-\Delta_{\mathbb{S}^{N-1}}$. It is well known that \begin{equation}\label{deflklk}
    \lambda_k=k(N-2+k),\quad k=0,1,2,\ldots,
    \end{equation}
    whose multiplicity is
    \[
    l_k:=\frac{(N+2k-2)(N+k-3)!}{(N-2)!k!}
    \]
    (note that $l_0=1$) and that \[
    \mathrm{Ker}(\Delta_{\mathbb{S}^{N-1}}+\lambda_k)
    =\mathbb{Y}_k(\mathbb{R}^N)|_{\mathbb{S}^{N-1}},
    \]
    where $\mathbb{Y}_k(\mathbb{R}^N)$ is the space of all homogeneous harmonic polynomials of degree $k$ in $\mathbb{R}^N$. It is standard that $\lambda_0=0$ and the corresponding eigenfunction of (\ref{deflk}) is the constant function that is $\Psi_{0,1}=c\in\mathbb{R}\setminus\{0\}$.

    Then the function $v$ is a solution of (\ref{rwlp}) if and only if $\varphi_{k,m}\in\mathcal{W}$ is a classical solution of the system
    \begin{eqnarray}\label{Ppwhl2p2tpy}
    \left\{ \arraycolsep=1.5pt
       \begin{array}{ll}
        (p-1)\varphi''_{k,m}
        +\left[(N-\mu-1)
        +\frac{(p-2)(N-s)}{1+\rho^{\frac{p-s+\mu}{p-1}}}\right]
        \frac{\varphi'_{k,m}}{\rho}
    -\frac{\lambda_k}{\rho^2}\varphi_{k,m} \\[4mm]
    + (r-1)C_{N,p,\mu,s}^{r-p}\left(\frac{N-p-\mu}{p-1}\right)^{2-p}
    \frac{\rho^{\frac{p-s+\mu}{p-1}-2}}
    {\left(1+\rho^{\frac{p-s+\mu}{p-1}}\right)^2}\varphi_{k,m}=0 \quad \mbox{in}\quad \rho\in(0,\infty),\\[4mm]
        \varphi'_{k,m}(0)=0 \quad\mbox{if}\quad k=0,\quad \mbox{and}\quad \varphi_{k,m}(0)=0 \quad\mbox{if}\quad k\geq 1,
        \end{array}
    \right.
    \end{eqnarray}
    for all $m=1,\ldots,l_k$, where $\mathcal{W}:=\{w\in C^1_{c,0}([0,\infty))| \int^\infty_0 |U'|^{p-2}|w'|^2 \rho^{N-1-\mu} \mathrm{d}\rho<\infty\}$.
    We make the change of variable $\rho=\tau^\varsigma$ where $\varsigma:=p/(p-s+\mu)>0$, and let
    \begin{equation}\label{Ppwhl2p2txyp}
    \eta_{k,m}(\tau)=\varphi_{k,m}(\rho),
    \end{equation}
    that transforms (\ref{Ppwhl2p2tpy}) into the following equations for all $\eta_{k,m}\in \widetilde{\mathcal{W}}$, $k=0,1,2,\ldots$ and $m=1,\ldots,l_k$,
    \begin{small}\begin{equation}\label{Ppwhl2p2tp}
    \eta''_{k,m}+\left(\frac{K-1}{p-1}
    +\frac{(p-2)K}{(p-1)(1+\tau^{\frac{p}{p-1}})}\right)
    \frac{\eta'_{k,m}}{\tau}
    -\frac{\varsigma^2\lambda_k}{(p-1)\tau^2}\eta_{k,m}
    +\frac{K(Kp-K+p)}{(p-1)^2}\frac{\tau^{\frac{p}{p-1}-2}}
    {(1+\tau^{\frac{p}{p-1}})^2}\eta_{k,m}=0.
    \end{equation}\end{small}
    where $\widetilde{\mathcal{W}}:=\{w\in C^1_{c,0}([0,\infty))| \int^\infty_0 |V'|^{p-2}|w'|^2 \tau^{K-1} \mathrm{d}\tau<\infty\}$, here $V(\tau)=U(\rho)$ and
    \begin{equation}
    K=\frac{p(N-s)}{p-s+\mu}>p.
    \end{equation}
    Here we have used the fact \begin{equation*}
    \varsigma^2(r-1)C_{N,p,\mu,s}^{r-p}\left(\frac{N-p-\mu}{p-1}\right)^{2-p}
    =\frac{K(Kp-K+p)}{p-1}.
    \end{equation*}

    Now, let us consider the linear operator
    \begin{equation}\label{lrntrpkb}
    \mathcal{A}_k(\eta):=
    \left(\tau^{K-1}|V'|^{p-2}\eta'\right)'
    +\frac{\tilde{p}^*-1}{p-1}\tau^{K-1}V^{\tilde{p}^*-2}
    \eta
    -\frac{\varsigma^2\lambda_k}{p-1} \tau^{K-3}|V'|^{p-2}\eta,\quad \eta\in \widetilde{\mathcal{W}}.
    \end{equation}
    Here $\tilde{p}^*=\frac{Kp}{K-p}$. Note that solving the equation \eqref{Ppwhl2p2tp} is equivalent to solve $\mathcal{A}_k(\eta)=0$ for all $k\geq 0$.

    $\bullet$ {\em The case $k=0$}.

    We know that the function
    \begin{equation*}
    \eta_0(\tau)=\frac{(p-1)-\tau^{\frac{p}{p-1}}}
    {(1+\tau^{\frac{p}{p-1}})^\frac{K}{p}}
    \end{equation*}
    solves the equation \eqref{Ppwhl2p2tp}. We claim that all the solutions are given by $\eta=c\eta_0$, $c\in\mathbb{R}$. Indeed, for $k=0$ we have that $\lambda_k=0$ and a straightforward computation shows that $\eta_0\in \widetilde{\mathcal{W}}$ and $\mathcal{A}_0(\eta_0)=0$. We look for a second linearly independent solution of the form
    \[
    w(\tau)=c(\tau)\eta_0(\tau).
    \]
    Then we get
    \[
    c''(r)\eta_0(\tau)+c'(r)\left[2\eta_0'
    +\frac{\eta_0}{\tau}\left(\frac{K-1}{p-1}
    +\frac{(p-2)K}{p-1}\frac{1}{1+\tau^{\frac{p}{p-1}}}\right)\right]
    =0,
    \]
    and hence
    \[
    \frac{c''(\tau)}{c'(\tau)}
    =-2\frac{\eta_0'(\tau)}{\eta_0(\tau)}
    -\frac{1}{\tau}\left(\frac{K-1}{p-1}
    +\frac{(p-2)K}{p-1}\frac{1}{1+\tau^{\frac{p}{p-1}}}\right).
    \]
    A direct computation shows that
    \[
    c'(\tau)=B\frac{\left(1+\tau^{\frac{p}{p-1}}\right)^{\frac{K(p-2)}{p}}}
    {(\eta_0(\tau))^2\tau^{\frac{K(p-2)+K-1}{p-1}}},
    \quad \mbox{for some}\quad B\in\mathbb{R}\setminus\{0\}.
    \]
    Therefore,
    \[c(\tau)\sim B\tau^{\frac{K-p}{p-1}}\quad \mbox{and}\quad w(\tau)=c(\tau)\eta_0(\tau)\sim B\quad \mbox{as}\quad \tau\to +\infty.\]
    However, $w\notin \widetilde{\mathcal{W}}$ by using the facts in \cite[Lemma 2.1]{PV21}.

$\bullet$ {\em The case $k\geq 1$}.

    In this case, we claim that all the solutions in $\widetilde{\mathcal{W}}$ of $\mathcal{A}_k(\eta)=0$ are identically zero. Assume there exists a function $\eta_{k}\in \widetilde{\mathcal{W}}$ such that $\mathcal{A}_k(\eta_{k})=0$, i.e. for any $\tau\geq 0$
    \begin{equation}\label{lrntrpk20}
    \begin{split}
    \left(\tau^{K-1}|V'|^{p-2}\eta_{k}'\right)'
    +\frac{\tilde{p}^*-1}{p-1}\tau^{K-1}V^{\tilde{p}^*-2}
    \eta_{k}
    -\frac{\varsigma^2 \lambda_k}{p-1} \tau^{K-3}|V'|^{p-2}\eta_{k}=0.
    \end{split}
    \end{equation}
    We claim that $\eta_{k}\equiv 0$ if $k\geq 1$. We argue by contradiction. Without loss of generality, we suppose that there exists $\tau_k>0$ (possibly $+\infty$) such that $\eta_{k}(\tau)>0$ for any $\tau\in (0,\tau_k)$ and $\eta_{k}(\tau_k)=0$. In particular, $\eta_{k}'(\tau_k)\leq 0$. Now, let $\eta_{*}\in \widetilde{\mathcal{W}}\setminus\{0\}$ be a solution of the following equation
    \begin{equation}\label{lrntrpk21}
    \begin{split}
    \left(\tau^{K-1}|V'|^{p-2}\eta_{*}'\right)'
    +\frac{\tilde{p}^*-1}{p-1}\tau^{K-1}V^{\tilde{p}^*-2}
    \eta_{*}
    -\frac{K-1}{p-1} \tau^{K-3}|V'|^{p-2}\eta_{*}=0,\quad \mbox{for all}\ \tau>0.
    \end{split}
    \end{equation}
    From \cite{PV21} we know equation \eqref{lrntrpk21} in $\widetilde{\mathcal{W}}$ admits only one solution given by $cV'$, so we take $\eta_*=V'$.
    Multiplying (\ref{lrntrpk20}) by $\eta_{*}$, (\ref{lrntrpk21}) by $\eta_{k}$, and integrating between $0$ and $\tau_k$ then subtracting the two expressions, we obtain
    \begin{align}\label{lrntrpk2s}
    \frac{\varsigma^2\lambda_k-(K-1)}{p-1}
    \int^{\tau_k}_0\tau^{K-3}|V'|^{p-2}\eta_{k} \eta_{*}\mathrm{d}\tau = & \int^{\tau_k}_0\left(\tau^{K-1}|V'|^{p-2}\eta_{k}'\right)' \eta_{*}\mathrm{d}\tau
    \nonumber\\
    & - \int^{\tau_k}_0\left(\tau^{K-1}|V'|^{p-2}\eta_{*}'\right)' \eta_{k}\mathrm{d}\tau
    \nonumber\\
    = & \tau_k^{K-1}|V'(\tau_k)|^{p-2}\eta_{k}'(\tau_k)\eta_{*}(\tau_k)\geq 0.
    \end{align}
    Note that under the assumption \eqref{cknib1c}, it is easy to verify that
    \begin{equation*}
    0<\left(\frac{p-s+\mu}{p}\right)^2\left[\frac{p(N-s)}{p-s+\mu}-1\right]
    <N-1,
    \end{equation*}
    thus $\varsigma^2\lambda_k>K-1$ for all $k\geq 1$.
    Then a contradiction arises in \eqref{lrntrpk2s} since $\eta_{k}'(\tau_k)\leq 0$, and $\eta_*(\tau)<0$ for any $\tau>0$. Here, we integrate by part and we use that $\eta_{k}(\tau_k)=0$. Thus all the solutions in $\widetilde{\mathcal{W}}$ of $\mathcal{A}_k(\eta_{k})=0$ are $\eta_{k}\equiv 0$ for $k\geq 1$.

    To sum up, we deduce that \eqref{Ppwhl2p2tpy} only admits one solution $c\varphi_0$ with $c\in \mathbb{R}$, where
    \begin{equation}\label{Ppwhl2pyfp}
    \varphi_0(\rho)=\frac{(p-1)-\rho^{\frac{p-s+\mu}{p-1}}}
    {(1+\rho^{\frac{p-s+\mu}{p-1}})^\frac{N-s}{p-s+\mu}}.
    \end{equation}
    That is, the space of solutions of (\ref{rwlp}) has dimension $1$ and is spanned by
    \begin{equation*}
    W_0(x)=\frac{(p-1)-|x|^{\frac{p-s+\mu}{p-1}}}
    {(1+|x|^{\frac{p-s+\mu}{p-1}})^\frac{N-s}{p-s+\mu}}.
    \end{equation*}
    Since $W_0\sim \frac{\partial U_\lambda}{\partial \lambda}|_{\lambda=1}=\frac{N-p-\mu}{p}U+x\cdot \nabla U$, the proof of of Theorem \ref{coroPpwhlpa} is complete.
    \qed

    \section{{\bfseries Spectral analysis}}\label{sectspana}

    Let us consider the following eigenvalue problem
    \begin{equation}\label{pevp}
    \begin{split}
    & \mathcal{L}_{U} [v]=\xi|x|^{-s} U^{r-2}v \quad \mbox{in}\quad \mathbb{R}^N,\quad v\in \mathcal{D}^{1,2}_{\mu,*}(\mathbb{R}^N),
    \end{split}
    \end{equation}
    where
    \begin{equation*}
    \begin{split}
    & \mathcal{L}_{U} [v]:=-\mathrm{div}(|x|^{-\mu}|\nabla U|^{p-2}\nabla v)-(p-2)\mathrm{div}(|x|^{-\mu}|\nabla U|^{p-4}(\nabla U\cdot\nabla v)\nabla U).
    \end{split}
    \end{equation*}
    As shown in this Proposition \ref{propcet}, this operator $\mathcal{L}_{U}$ has a discrete spectrum. By the standard Sturm-Liouville theory as in \cite{FN19,FZ22}, and from Theorem \ref{coroPpwhlpa} we can classify the eigenspaces corresponding to $\mathcal{L}_{U}$. It is easy to verify that the functions $U$ and $(\frac{N-p-\mu}{p}U
    +x\cdot \nabla U)$ corresponding to $\xi=(p-1)$ and $(r-1)$, furthermore, since $U$ is positive and $(\frac{N-p-\mu}{p}U
    +x\cdot \nabla U)$ has only one zero, the classical Sturm-Liouville theory ensures that $(p-1)$ and $(r-1)$ are the first two eigenvalues, with multiplicity $1$ due to (CKN) inequality \eqref{cknp} and Theorem \ref{coroPpwhlpa}.

    Indeed, as in Section \ref{sectndr}, given an eigenfunction of the form $v(x)=\varphi(\rho)\Psi(\theta)$ where $\rho=|x|$ and $\theta=\frac{x}{|x|}\in \mathbb{S}^{N-1}$, the eigenvalue problem corresponds to the following system
    \begin{align}\label{epcse}
    0= & \Delta_{\mathbb{S}^{N-1}}\varphi+\lambda  \varphi\quad \mbox{on}\ \mathbb{S}^{N-1}, \\ \label{epcs}
    0= & (p-1)\varphi''
        +\left[(N-\mu-1)
        +\frac{(p-2)(N-s)}{1+\rho^{\frac{p-s+\mu}{p-1}}}\right]
        \frac{\varphi'}{\rho}
    -\frac{\lambda}{\rho^2}\varphi \nonumber\\
    & + \xi C_{N,p,\mu,s}^{r-p}\left(\frac{N-p-\mu}{p-1}\right)^{2-p}
    \frac{\rho^{\frac{p-s+\mu}{p-1}-2}}
    {\left(1+\rho^{\frac{p-s+\mu}{p-1}}\right)^2}\varphi\quad \mbox{on}\ \rho\in[0,\infty).
    \end{align}
    Making the change of variable: $\rho=\tau^\varsigma$ with $\varsigma:=p/(p-s+\mu)$,  and let
    \begin{equation*}
    \eta(\tau)=\varphi(\rho),
    \end{equation*}
    then multiplying by the integrating factor $\tau^{K-1}$, \eqref{epcs} is equivalent to
    \begin{align}\label{epcsb}
    0= \left((p-1)|V'|^{p-2}\tau^{K-1}\eta'\right)'
    +\varsigma^2\xi V^{\tilde{p}^*-2}\tau^{K-1}\eta
    - \varsigma^2\lambda |V'|^{p-2}\tau^{K-3}\eta \quad \mbox{on}\ \tau\in[0,\infty),
    \end{align}
    where $V(\tau)=U(\rho)$, $K=\frac{p(N-s)}{p-s+\mu}>p$ satisfying $r=\frac{Kp}{K-p}=\tilde{p}^*$. For each $\lambda$ the ordinary differential equation \eqref{epcsb} takes the form of the Sturm-Liouville eigenvalue problem
    \begin{align}\label{epcsbsl}
    L\eta+ \alpha\eta=0\quad \mbox{on}\ [0,\infty),
    \end{align}
    where $\alpha=\varsigma^{2}\xi$, and
    \begin{align*}
    L\eta=\frac{1}{\mathfrak{W}}[(\mathfrak{P}\eta')'
    -\mathfrak{Q}\eta]
    \end{align*}
    with
    \begin{align*}
    \mathfrak{P}(\tau)=(p-1)|V'|^{p-2}\tau^{K-1},\quad
    \mathfrak{Q}(\tau)=\varsigma^2\lambda |V'|^{p-2}\tau^{K-3},\quad
    \mathfrak{W}(\tau)=V^{\tilde{p}^*-2}\tau^{K-1},
    \end{align*}
    and the eigenfunctions belong to
    \begin{align*}
    \mathcal{H}:=\{g: [0,\infty)\mapsto \mathbb{R}: g\in L^2([0,\infty);\mathfrak{W}), g'\in L^2([0,\infty);\mathfrak{P})\}.
    \end{align*}
    When $K$ is an integer, \cite[Lemma B.3]{FZ22} states that,
    \begin{itemize}
    \item[$(1)$]
    if $\eta_1$ and $\eta_2$ are two eigenfunctions corresponding to the same eigenvalue $\alpha$, then $\eta_1=c\eta_2$ for some $c\in\mathbb{R}$;
    \item[$(2)$]
    the $i$-th eigenfunction of $L$ has $i-1$ interior zeros.
    \end{itemize}
    On the other hand, since \eqref{epcsbsl} is an ODE, even if $K$ is not an integer the conclusion also holds.
    Therefore we can deduce directly the following conclusion.

    \begin{proposition}\label{propev}
    Assume that \eqref{cknib1c} holds. Let $\xi_i$, $i=1,2,\ldots,$ denote the eigenvalues of (\ref{pevp}) in increasing order. Then $\xi_1=(p-1)$ is simple and the corresponding eigenfunction is spanned by $U$, $\xi_2=r-1$ and the corresponding eigenfunction is spanned by $(\frac{N-p-\mu}{p}U+x\cdot\nabla U)$.
    \end{proposition}

    \begin{proof}
    From the statements as previous in this section, we have $U$ is an eigenfunction of $\xi_1=p-1$ and also $(\frac{N-p-\mu}{p}U+x\cdot\nabla U)$ is an eigenfunction of $\xi_2=r-1$.

    Same as in the proof of Lemma \ref{lemcztj}, we assume by approximation that $\varphi\in C^1_{c,0}(\mathbb{R}^N)$ (since $C^1_{c,0}(\mathbb{R}^N)$ is dense in $\mathcal{D}^{1,2}_{\mu,*}(\mathbb{R}^N)$) and define
    \[
    \mathfrak{F}(u):=\int_{\mathbb{R}^N}|x|^{-\mu}|\nabla u|^p \mathrm{d}x
    -\mathcal{S}\left(\int_{\mathbb{R}^N}
    |x|^{-s}|u|^{r} \mathrm{d}x\right)^\frac{p}{r}.
    \]
    We know $U$ is a local minimum of the functional $\mathfrak{F}$, then
    \begin{align*}
    0\leq & \frac{\mathrm{d}^2}{\mathrm{d}\varepsilon^2}\Big|_{\varepsilon=0} \mathfrak{F}(U+\varepsilon\varphi)
    \\ = & p \int_{\mathbb{R}^N}|x|^{-\mu}|\nabla U|^{p-2}|\nabla \varphi|^2\mathrm{d}x
    + p(p-2) \int_{\mathbb{R}^N}|x|^{-\mu}|\nabla U|^{p-4}(\nabla U\cdot\nabla\varphi)^2\mathrm{d}x
    \\ & - \mathcal{S}
    \Bigg[
    p\left(p-r\right)
    \left(\int_{\mathbb{R}^N}|x|^{-s}|U|^{r} \mathrm{d}x\right)^{\frac{p}{r}-2}
    \left(\int_{\mathbb{R}^N}
    |x|^{-s}|U|^{r-2}U\varphi\mathrm{d}x\right)^2
    \\ & \quad\quad + p(r-1)
    \left(\int_{\mathbb{R}^N}|x|^{-s}|U|^{r} \mathrm{d}x\right)^{\frac{p}{r}-1}
    \int_{\mathbb{R}^N}
    |x|^{-s}|U|^{r-2}\varphi^2\mathrm{d}x
    \Bigg].
    \end{align*}
    Noting that
    \begin{align}\label{xi1a}
    \left(\int_{\mathbb{R}^N}
    |x|^{-s}|U|^{r-2}U\varphi\mathrm{d}x\right)^2
    \leq \int_{\mathbb{R}^N}
    |x|^{-s}|U|^{r}\mathrm{d}x
    \int_{\mathbb{R}^N}
    |x|^{-s}|U|^{r-2}\varphi^2\mathrm{d}x,
    \end{align}
    which implies that
    \begin{align}\label{xi1b}
    \left(p-1\right)
    \int_{\mathbb{R}^N}|x|^{-s}|U|^{r-2}\varphi^2\mathrm{d}x
    \leq & \int_{\mathbb{R}^N}|x|^{-\mu}|\nabla U|^{p-2}|\nabla \varphi|^2\mathrm{d}x
    \nonumber\\
    & + (p-2) \int_{\mathbb{R}^N}|x|^{-\mu}|\nabla U|^{p-4}(\nabla U\cdot\nabla\varphi)^2\mathrm{d}x,
    \end{align}
    due to
    \[
    \int_{\mathbb{R}^N}|x|^{-s}|U|^{r} \mathrm{d}x
    =\mathcal{S}^{\frac{r}{r-p}}.
    \]
    The equality achieves in \eqref{xi1a} if and only if $\varphi=cU$ which implies also as in \eqref{xi1b}, therefore the first eigenvalue $\xi_1=p-1$ with corresponding eigenfunction $U$ (up to scalar multiplications). Furthermore, Theorem \ref{coroPpwhlpa} indicates $(\frac{N-p-\mu}{p}U+x\cdot\nabla U)$ is the unique  eigenfunction (up to scalar multiplications) of $\xi_2=r-1$. Now the proof is completed.
    \end{proof}

    In particular, Proposition \ref{propev} implies
    \begin{align}\label{czkj}
    T_{U}\mathcal{M}=\mathrm{Span}\left\{U, \frac{N-p-\mu}{p}U+x\cdot\nabla U\right\},
    \end{align}
    where $\mathcal{M}=\{cU_{\lambda}: c\in\mathbb{R}, \lambda>0\}$ is the set of extremal functions for (CKN) inequality \eqref{cknp}.
    From Proposition \ref{propev}, we directly obtain
    \begin{proposition}\label{propevl}
    Assume that \eqref{cknib1c} holds. Then there exists $\tau=\tau(N,p,\mu,s)>0$ such that for any function $\mathcal{D}^{1,2}_{\mu,*}(\mathbb{R}^N)$ orthogonal to $T_{U} \mathcal{M}$, it holds that
    \begin{small}\begin{align*}
    \int_{\mathbb{R}^N}|x|^{-\mu}\left[|\nabla U|^{p-2}|\nabla v|^2+(p-2)|\nabla U|^{p-4}(\nabla U\cdot\nabla v)^2\right]\mathrm{d}x
    \geq \left[(r-1)+2\tau\right]
    \int_{\mathbb{R}^N}|x|^{-s}U^{r-2}v^2\mathrm{d}x.
    \end{align*}\end{small}
    \end{proposition}

    Following \cite{FZ22}, we give the following remark which will  be important to give a meaning to the notion of  ``orthogonal to $T_{U} \mathcal{M}$" for functions which are not necessarily in $\mathcal{D}^{1,2}_{\mu,*}(\mathbb{R}^N)$. Note that Proposition \ref{propcet} indicates $\mathcal{D}^{1,2}_{\mu,*}(\mathbb{R}^N)\subset L^2_{s,*}(\mathbb{R}^N)$.

    \begin{remark}\label{remd12p}
    For any $w\in T_{U} \mathcal{M}$ it holds $U^{r-2}w \in L_s^{\frac{r}{r-1}}(\mathbb{R}^N)
    =\left(L_s^{r}(\mathbb{R}^N)\right)'$, here $L_s^{q}(\mathbb{R}^N)$ is the set of measurable functions with the norm $\|\varphi\|_{L_s^{q}(\mathbb{R}^N)}:
    =\left(\int_{\mathbb{R}^N}|x|^{-s}|\varphi|^{q} \mathrm{d}x\right)^{\frac{1}{q}}$.
    Hence, by abuse of notation, for any function $v\in L_s^{r}(\mathbb{R}^N)$ we say that $v$ is orthogonal to $T_{U} \mathcal{M}$ in $L^2_{s,*}(\mathbb{R}^N)$ if \[
    \int_{\mathbb{R}^N}|x|^{-s}U^{r-2} w v\mathrm{d}x=0,\quad \mbox{for all}\quad w \in T_{U} \mathcal{M},
    \]
    in particular, taking $w=U$, we deduce
    \[
    \int_{\mathbb{R}^N}|x|^{-\mu} |\nabla U|^{p-2}\nabla U\cdot \nabla v \mathrm{d}x=\int_{\mathbb{R}^N}|x|^{-s}U^{r-1} v\mathrm{d}x=0.
    \]

    \end{remark}


    Note that, by H\"{o}lder inequality and (CNK) inequality \eqref{cknp}, $\mathcal{D}^{1,p}_{\mu}(\mathbb{R}^N)\subset
    L_s^{r}(\mathbb{R}^N)\subset L^2_{s,*}(\mathbb{R}^N)$ if $r\geq 2$.
    Hence, the notion of orthogonality introduced above is particularly relevant when $r<2$ (equivalently, $p<\frac{2(N-\mu)}{N+2-s}$).

    \subsection{{\bfseries Improvements of spectral gap}}\label{subsectsge}

    Same as \cite[Proposition 3.8]{FZ22}, the goal of this subsection is improving the spectral gap obtained in  Proposition \ref{propevl}, that is, we will give the following spectral gap-type estimates for $p\geq 2$ and $1<p<2$, respectively.

    \begin{lemma}\label{lemsgap}
    Assume that \eqref{cknib1c} holds and $p\geq 2$. Given any $\gamma_0>0$, there exists $\overline{\delta}=\overline{\delta}(N,p,\mu,s,\gamma_0)>0$ such that for any function $v\in \mathcal{D}^{1,p}_{\mu}(\mathbb{R}^N)$ orthogonal to $T_{U} \mathcal{M}$ in $L^2_{s,*}(\mathbb{R}^N)$ satisfying $\|v\|_{\mathcal{D}^{1,p}_\mu(\mathbb{R}^N)}\leq \overline{\delta}$, we have
    \begin{equation*}
    \begin{split}
    & \int_{\mathbb{R}^N}|x|^{-\mu}\left[
    |\nabla U|^{p-2}|\nabla v|^2+(p-2)|\bar{\omega}|^{p-2}(|\nabla (U+v)|-|\nabla U|)^2\right]\mathrm{d}x \\
    & \quad \geq \left[(r-1)+\tau\right]
    \int_{\mathbb{R}^N}|x|^{-s}U^{r-2}|v|^2\mathrm{d}x,
    \end{split}
    \end{equation*}
    where $\tau>0$ is given in Proposition \ref{propevl}, and $\bar{\omega}: \mathbb{R}^{2N}\to \mathbb{R}^N$ is defined in analogy to Lemma \ref{lemui1p}:
    \begin{eqnarray*}
    \bar{\omega}=\bar{\omega}(\nabla U,\nabla (U+v))=
    \left\{ \arraycolsep=1.5pt
       \begin{array}{ll}
        \nabla U,\ \ &{\rm if}\ \ |\nabla U|<|\nabla (U+v)|,\\[3mm]
        \left(\frac{|\nabla (U+v)|}{|\nabla U|}\right)^{\frac{1}{p-2}}\nabla (U+v),\ \ &{\rm if}\ \  |\nabla (U+v)|\leq |\nabla U|.
        \end{array}
    \right.
    \end{eqnarray*}
    \end{lemma}

    \begin{proof}
    The proof can be found in \cite[Lemma 4.6]{DT22}, so we omit it. 
    \end{proof}

    \begin{lemma}\label{lemsgap2}
    Assume that \eqref{cknib1c} holds and $1<p<2$. Given any $\gamma_0>0$, $\mathcal{C}>0$ there exists $\overline{\delta}=\overline{\delta}(N,p,\mu,s,\gamma_0,\mathcal{C})>0$ such that for any function $v\in \mathcal{D}^{1,p}_{\mu}(\mathbb{R}^N)$ orthogonal to $T_{U} \mathcal{M}$ in $L^2_{s,*}(\mathbb{R}^N)$ satisfying $\|v\|_{\mathcal{D}^{1,p}_{\mu}(\mathbb{R}^N)}\leq \overline{\delta}$, the following holds:
    \begin{itemize}
    \item[$(i)$]
    when $1<p\leq\frac{2(N-\mu)}{N+2-s}$, we have
    \begin{small}
    \begin{align*}
    & \int_{\mathbb{R}^N}
    |x|^{-\mu}\left[
    |\nabla U|^{p-2}|\nabla v|^2
    +(p-2)|\tilde{\omega}|^{p-2}(|\nabla (U+v)|-|\nabla U|)^2
    +\gamma_0 \min\{|\nabla v|^p,|\nabla U|^{p-2}|\nabla v|^2\}
    \right]\mathrm{d}x \\
    & \quad \geq \left[(r-1)+\tau\right]
    \int_{\mathbb{R}^N}|x|^{-s}
    \frac{(U+\mathcal{C}|v|)^{r}}{U^2+|v|^2}|v|^2\mathrm{d}x;
    \end{align*}
    \end{small}
    \item[$(ii)$]
    when $\frac{2(N-\mu)}{N+2-s}< p<2$, we have
    \begin{small}
    \begin{align*}
    & \int_{\mathbb{R}^N}
    |x|^{-\mu}\left[
    |\nabla U|^{p-2}|\nabla v|^2
    +(p-2)|\tilde{\omega}|^{p-2}(|\nabla (U+v)|-|\nabla U|)^2
    +\gamma_0 \min\{|\nabla v|^p,|\nabla U|^{p-2}|\nabla v|^2\}
    \right]\mathrm{d}x \\
    & \quad \geq \left[(r-1)+\tau\right]
    \int_{\mathbb{R}^N}|x|^{-s}U^{r-2}|v|^2\mathrm{d}x,
    \end{align*}
    \end{small}
    \end{itemize}
    where $\tau>0$ is given in Proposition \ref{propevl}, and $\tilde{\omega}: \mathbb{R}^{2N}\to \mathbb{R}^N$ is defined in analogy to Lemma \ref{lemui1p}:
    \begin{eqnarray*}
    \tilde{\omega}=\tilde{\omega}(\nabla U,\nabla (U+v))=
    \left\{ \arraycolsep=1.5pt
       \begin{array}{ll}
        \left(\frac{|\nabla (U+v)|}{(2-p)|\nabla (U+v)|+(p-1)|\nabla U|}\right)^{\frac{1}{p-2}}\nabla U,\ \ &{\rm if}\ \  |\nabla U|<|\nabla (U+v)|\\[3mm]
        \nabla U,\ \ &{\rm if}\ \ |\nabla (U+v)|\leq |\nabla U|
        \end{array}.
    \right.
    \end{eqnarray*}
    \end{lemma}

    \begin{proof}
    The proof is similar to the proof of Lemma \ref{lemsgap}, but it is more complicated. Furthermore, this proof can be deduced with minor changes as in \cite[Lemma 4.9]{DT22} but it deals with the radial case, so here we will give the proof directly by using the compact embedding as in Lemma \ref{propcetl}. We argue by contradiction in these two cases.

    $\bullet$ {\em The case $1<p\leq\frac{2(N-\mu)}{N+2-s}$} which implies $r\leq2$.
    Suppose the inequality does not hold, then there exists a sequence $0\not\equiv v_n\to 0$ in $\mathcal{D}^{1,p}_{\mu}(\mathbb{R}^N)$, with $v_n$ orthogonal to $T_{U} \mathcal{M}$ in $L^2_{s,*}(\mathbb{R}^N)$, such that
    \begin{small}
    \begin{align}\label{evbc2}
    & \int_{\mathbb{R}^N}|x|^{-\mu}\Big[
    |\nabla U|^{p-2}|\nabla v_n|^2
    +(p-2)|\tilde{\omega}_n|^{p-2}(|\nabla (U+v_n)|-|\nabla U|)^2
    \nonumber\\& \quad \quad +\gamma_0 \min\{|\nabla v_n|^p,|\nabla U|^{p-2}|\nabla v_n|^2\}
    \Big]\mathrm{d}x \nonumber\\
    & \quad <\left[(r-1)+\tau\right]
    \int_{\mathbb{R}^N}|x|^{-s}
    \frac{(U+\mathcal{C}|v_n|)^{r}}{U^2+|v_n|^2}|v_n|^2\mathrm{d}x,
    \end{align}
    \end{small}
    where $\omega_n$ corresponds to $v_n$ as in the statement.
    Let
    \begin{align}\label{defevh}
    \varepsilon_n:=\left(\int_{\mathbb{R}^N}
    |x|^{-\mu}(|\nabla U|+|\nabla v_n|)^{p-2}|\nabla v_n|^2 \mathrm{d}x\right)^{\frac{1}{2}},\quad \widehat{v}_n=\frac{v_n}{\varepsilon_n}.
    \end{align}
    Note that, since $1<p<2$, it follows by H\"{o}lder inequality that
    \begin{small}\begin{align*}
    \int_{\mathbb{R}^N}
    |x|^{-\mu}(|\nabla U|+|\nabla v_n|)^{p-2}|\nabla v_n|^2 dx\leq & \int_{\mathbb{R}^N}|x|^{-\mu}|\nabla v_n|^{p-2}|\nabla v_n|^2\mathrm{d}x
    = \int_{\mathbb{R}^N}|x|^{-\mu}|\nabla v_n|^{p}\mathrm{d}x
    \to 0,
    \end{align*}\end{small}
    hence $\varepsilon_n\to 0$, as $n\to \infty$.

    Since the integrand in the left hand side of \eqref{evbc2} is nonnegative, for any $R>1$ we have
    \begin{align}\label{evbcb2}
    & \int_{B_R\backslash B_{\frac{1}{R}}}
    |x|^{-\mu}\bigg[
    |\nabla U|^{p-2}|\nabla \widehat{v}_n|^2
    +(p-2)|\tilde{\omega}_n|^{p-2}\left(\frac{|\nabla (U+v_n)|-|\nabla U|}{\varepsilon_n}\right)^2 \nonumber\\
    & \quad\quad + \gamma_0 \min\{\varepsilon_n^{p-2}|\nabla \widehat{v}_n|^p,|\nabla U|^{p-2}|\nabla \widehat{v}_n|^2\}
    \bigg]\mathrm{d}x \nonumber\\
    & \quad < \left[(r-1)+\tau\right]
    \int_{\mathbb{R}^N}|x|^{-s}\frac{(U+\mathcal{C}|v_n|)^{r}}
    {U^2+|v_n|^2}|\widehat{v}_n|^2\mathrm{d}x.
    \end{align}
    Now, same as the proof of Lemma \ref{lemsgap}, let us
    fix $R>1$ which can be chosen arbitrarily large, and set
    \begin{align}\label{defrsi}
    \mathcal{R}_n:=\{2|\nabla U|\geq |\nabla v_n|\}&,\quad \mathcal{S}_n:=\{2|\nabla U|< |\nabla v_n|\},  \nonumber\\
    \mathcal{R}_{n,R}:=\left(B_R\backslash B_{\frac{1}{R}}\right)\cap \mathcal{R}_n&,\quad
    \mathcal{S}_{n,R}:=\left(B_R\backslash B_{\frac{1}{R}}\right)\cap \mathcal{S}_n,
    \end{align}
    thus $B_R\backslash B_{\frac{1}{R}}=\mathcal{R}_{n,R} \cup \mathcal{S}_{n,R}$.
    From \cite[(2.2)]{FZ22}, that is, for $1<p<2$, there exists $c(p)>0$ such that
    \begin{equation}\label{evbcb2i}
    p|x|^{p-2}|y|^2+p(p-2)|\tilde{\omega}|^{p-2}(|x|-|x+y|)^2\geq c(p)\frac{|x|}{|x|+|y|}|x|^{p-2}|y|^2,
    \end{equation}
     for all $|x|\neq 0$ and $y\in\mathbb{R}^N$, then we have
    \begin{align*}
    & \left|\frac{\nabla U}{\varepsilon_n}\right|^{p-2}|\nabla \widehat{v}_n|^2
    +(p-2)\left|\frac{\tilde{\omega}_n}{\varepsilon_n}\right|^{p-2}
    \left(\left|\frac{\nabla U }{\varepsilon_n}+\nabla \widehat{v}_n\right|-\left|\frac{\nabla  U}{\varepsilon_n}\right|\right)^2 \\
    & \quad \geq c(p)\frac{|\nabla U|/\varepsilon_n}{|\nabla U|/\varepsilon_n+|\nabla \widehat{v}_n|}\left|\frac{\nabla  U}{\varepsilon_n}\right|^{p-2}|\nabla \widehat{v}_n|^2,
    \end{align*}
    then,
    \begin{small}\begin{align*}
    \left| \nabla U \right|^{p-2}|\nabla \widehat{v}_n|^2
    +(p-2)\left|\tilde{\omega}_n \right|^{p-2}\left(\frac{|\nabla (U+v_n)|-|\nabla  U| }{\varepsilon_n}\right)^2
     \geq c(p) \left| \nabla  U \right|^{p-2}|\nabla \widehat{v}_n|^2,\quad \mbox{in}\quad \mathcal{R}_{n,R}.
    \end{align*}\end{small}
    Therefore, combining this bound with \eqref{evbcb2}, we obtain
    \begin{align}\label{evbcb2l}
    & c(p)\int_{\mathcal{R}_{n,R}}
    |x|^{-\mu}\left| \nabla  U \right|^{p-2}|\nabla \widehat{v}_n|^2 \mathrm{d}x
    +\gamma_0\varepsilon^{p-2}_n\int_{\mathcal{S}_{n,R}}
    |x|^{-\mu}|\nabla \widehat{v}_n|^p \mathrm{d}x \nonumber\\
    \leq & \int_{B_R\backslash B_{\frac{1}{R}}}
    |x|^{-\mu}\bigg[
    |\nabla U|^{p-2}|\nabla \widehat{v}_n|^2
    +(p-2)|\tilde{\omega}_n|^{p-2}\left(\frac{|\nabla (U+v_n)|-|\nabla U|}{\varepsilon_n}\right)^2 \nonumber\\
    & \quad\quad + \gamma_0 \min\{\varepsilon_n^{p-2}|\nabla \widehat{v}_n|^p,|\nabla U|^{p-2}|\nabla \widehat{v}_n|^2\}
    \bigg]\mathrm{d}x \nonumber\\
    < & \left[(r-1)+\tau\right]
    \int_{\mathbb{R}^N}|x|^{-s}
    \frac{(U+\mathcal{C}|v_n|)^{r}}{U^2+|v_n|^2}|\widehat{v}_n|^2
    \mathrm{d}x.
    \end{align}
    In particular, this implies that
    \begin{align}\label{evbcb2li}
    1= & \varepsilon^{-2}_n\int_{\mathbb{R}^N}
    |x|^{-\mu}(|\nabla U|+|\nabla v_n|)^{p-2}|\nabla v_n|^2 dx \nonumber\\
    \leq & C(p)\left[\int_{\mathcal{R}_{n}}
    |x|^{-\mu}\left| \nabla  U \right|^{p-2}|\nabla \widehat{v}_n|^2 \mathrm{d}x
    +\varepsilon^{p-2}_n\int_{\mathcal{S}_{n}}
    |x|^{-\mu}|\nabla \widehat{v}_n|^p \mathrm{d}x\right] \nonumber\\
    \leq & C(N,p,\gamma_0)\left[(r-1)+\tau\right]
    \int_{\mathbb{R}^N}|x|^{-s}
    \frac{(U+\mathcal{C}|v_n|)^{r}}{U^2+|v_n|^2}|\widehat{v}_n|^2
    \mathrm{d}x.
    \end{align}
    Furthermore, thanks to \eqref{propcetlcpi} in Corollary \ref{propcetlpi}, for $n$ large enough so that $\varepsilon_n$ small we have
    \begin{align}\label{evbcb2lib}
    \int_{\mathbb{R}^N}|x|^{-s}
    \frac{(U+\mathcal{C}|v_n|)^{r}}{U^2+|v_n|^2}|\widehat{v}_n|^2
    \mathrm{d}x
    \leq & C(N,p,\mathcal{C})\int_{\mathbb{R}^N}|x|^{-s}
    (U+|v_n|)^{r-2}|\widehat{v}_n|^2\mathrm{d}x \nonumber\\
    \leq & C(N,p,\mu,s,\mathcal{C})\int_{\mathbb{R}^N}
    |x|^{-\mu}(|\nabla U|+|\nabla v_n|)^{p-2}|\nabla \widehat{v}_n|^2\mathrm{d}x \nonumber\\
    = & C(N,p,\mu,s,\mathcal{C}).
    \end{align}
    Hence, combining \eqref{evbcb2l} with \eqref{evbcb2lib}, by the definition of $\mathcal{S}_{i,R}$ we have
    \begin{small}\begin{align*}
    \varepsilon^{-2}_n\int_{\mathcal{S}_{n,R}}
    |x|^{-\mu}|\nabla U|^p \mathrm{d}x
    \leq \frac{\varepsilon^{-2}_n}{2^p}\int_{\mathcal{S}_{n,R}}
    |x|^{-\mu}|\nabla v_n|^p \mathrm{d}x
    = \frac{\varepsilon^{p-2}_n}{2^p}\int_{\mathcal{S}_{n,R}}
    |x|^{-\mu}|\nabla \widehat{v}_n|^p \mathrm{d}x
    \leq C,
    \end{align*}\end{small}
    then since
    \[
    0<c(R)\leq |\nabla U|\leq C(R)\quad \mbox{inside}\quad B_R\backslash B_{\frac{1}{R}},\quad \mbox{for all}\quad  R>1,
    \]
    for some constants $c(R)\leq C(R)$ depending only on $R$, we conclude that
    \begin{equation}\label{csirt0}
    |\mathcal{S}_{n,R}|\to 0\quad \mbox{as}\quad n\to \infty,\quad \mbox{for all}\quad R>1.
    \end{equation}
    Now, from \eqref{defevh} we have
    \[
    \int_{\mathbb{R}^N}
    |x|^{-\mu}(|\nabla U|+\varepsilon_n|\nabla \widehat{v}_n|)^{p-2}|\nabla \widehat{v}_n|^2 \mathrm{d}x\leq 1,
    \]
    then according to Lemma \ref{propcetl}, we deduce that $\widehat{v}_n\rightharpoonup \widehat{v}$ in $\mathcal{D}^{1,p}_\mu(\mathbb{R}^N)$ for some $\widehat{v}\in \mathcal{D}^{1,p}_{\mu}(\mathbb{R}^N)\cap L^2_{s,*}(\mathbb{R}^N)$, and
    \begin{equation}\label{csirtc}
    \int_{\mathbb{R}^N}|x|^{-s}
    \frac{(U+\mathcal{C}|v_n|)^{r}}{U^2+|v_n|^2}|\widehat{v}_n|^2
    \mathrm{d}x
    \to \int_{\mathbb{R}^N}|x|^{-s}U^{r-2}
    |\widehat{v}|^2\mathrm{d}x,
    \end{equation}
    as $n\to \infty$, for any $\mathcal{C}\geq 0$.
    Also, using \eqref{evbcb2l} and \eqref{evbcb2lib} again we have
    \begin{equation*}
    \int_{\mathcal{R}_{n,R}}
    |x|^{-\mu}\left| \nabla  U \right|^{p-2}|\nabla \widehat{v}_n|^2 \mathrm{d}x
    \leq C(N,p,\mu,s,\mathcal{C}),
    \end{equation*}
    therefore \eqref{csirt0} and $\widehat{v}_n\rightharpoonup \widehat{v}$ in $\mathcal{D}^{1,p}_\mu(\mathbb{R}^N)$ imply that, up to a subsequence,
    \[
    \widehat{v}_n \chi_{\mathcal{R}_{n,R}}\rightharpoonup \widehat{v} \chi_{B_R\backslash B_{\frac{1}{R}}}\quad \mbox{in}\quad \mathcal{D}^{1,2}_{\mu,*}(\mathbb{R}^N),\quad \mbox{for all}\quad R>1.
    \]
    Here $\chi_E$ denotes that $\chi_E=1$ if $x\in E$ and  $\chi_E=0$ if $x\notin E$.
    In addition, letting $n\to \infty$ in \eqref{evbcb2li} and \eqref{evbcb2lib}, and using \eqref{csirtc}, we deduce that
    \begin{equation}\label{csirtl1}
    0<c(N,p,\mu,s,\mathcal{C},\gamma_0)
    \leq \|\widehat{v}\|_{L^2_{s,*}(\mathbb{R}^N)}
    \leq C.
    \end{equation}
    Let us write
    \[
    \widehat{v}_n=\widehat{v}+\varphi_n,\quad\mbox{with}\quad \varphi_n:= \widehat{v}_n-\widehat{v},
    \]
    we have
    \[
    \varphi_n \rightharpoonup 0 \quad \mbox{in}\quad   \mathcal{D}^{1,p}_{\mu}(\mathbb{R}^N)
    \quad \mbox{and}\quad \varphi_n \chi_{\mathcal{R}_{n}}\rightharpoonup 0
    \quad \mbox{locally in} \quad \mathcal{D}^{1,2}_{\mu,*}(\mathbb{R}^N\setminus\{\mathbf{0}\}).
    \]
    We now look at the left side of \eqref{evbcb2}.
    The strong convergence $v_n\to 0$ in $\mathcal{D}^{1,p}_{\mu}(\mathbb{R}^N)$ implies that, $|\omega_n|\to |\nabla U|$ a.e. in $\mathbb{R}^N$. Then, let us rewrite
    \begin{align*}
    \left(\frac{|\nabla (U+v_n)|-|\nabla U|}{\varepsilon_n}\right)^2
    = & \left(\left[\int^1_0\frac{\nabla U+ t\nabla v_n}{|\nabla U+ t\nabla v_n|}\mathrm{d}t\right]\cdot \nabla \widehat{v}_n\right)^2 \\
    = & \left(\left[\int^1_0\frac{\nabla U+ t\nabla v_n}{|\nabla U+ t\nabla v_n|}\mathrm{d}t\right]\cdot \nabla (\widehat{v}+\varphi_n)\right)^2.
    \end{align*}
    Hence, if we set
    \[
    f_{n,1}=\left[\int^1_0\frac{\nabla U+ t\nabla v_n}{|\nabla U+ t\nabla v_n|}\mathrm{d}t\right]\cdot \nabla \widehat{v},\quad
    f_{n,2}=\left[\int^1_0\frac{\nabla U+ t\nabla v_n}{|\nabla U+ t\nabla v_n|}\mathrm{d}t\right]\cdot \nabla \varphi_n,
    \]
    since $\frac{\nabla U+ t\nabla v_n}{|\nabla U+ t\nabla v_n|}\to \frac{\nabla U}{|\nabla U|}$ a.e., it follows from Lebesgue's dominated convergence theorem that
    \[
    f_{n,1}\to \frac{\nabla U}{|\nabla U|}\cdot\nabla \widehat{v}\quad \mbox{locally in}\quad L^2(\mathbb{R}^N\backslash\{\mathbf{0}\}),\quad f_{n,2}\chi_{\mathcal{R}_n}\rightharpoonup 0\quad \mbox{locally in}\quad L^2(\mathbb{R}^N\backslash\{\mathbf{0}\}).
    \]
    Thus, the left hand side of \eqref{evbcb2} from below as follows:
    \begin{align}\label{evbcbb2}
    & \int_{\mathcal{R}_{n,R}}
    |x|^{-\mu}\left[|\nabla U|^{p-2}|\nabla \widehat{v}_n|^2+(p-2)|\tilde{\omega}_n|^{p-2}\left(\frac{|\nabla (U+v_n)|-|\nabla U|}{\varepsilon_n}\right)^2\right]\mathrm{d}x \nonumber\\
    = & \int_{\mathcal{R}_{n,R}}
    |x|^{-\mu}\left[|\nabla U|^{p-2}\left(|\nabla \widehat{v}|^2+2 \nabla \varphi_n\cdot \nabla \widehat{v}\right)+(p-2)|\tilde{\omega}_n|^{p-2}
    \left(f_{n,1}^2+2f_{n,1}f_{n,2}\right)\right]\mathrm{d}x \nonumber\\
    & + \int_{\mathcal{R}_{n,R}}
    |x|^{-\mu}\left[|\nabla U|^{p-2} |\nabla \varphi_n|^2+(p-2)|\tilde{\omega}_n|^{p-2}f_{n,2}^2\right]\mathrm{d}x
    \nonumber\\
    \geq & \int_{\mathcal{R}_{n,R}}
    |x|^{-\mu}\left[|\nabla U|^{p-2}\left(|\nabla \widehat{v}|^2+2 \nabla \varphi_n\cdot \nabla \widehat{v}\right)+(p-2)|\tilde{\omega}_n|^{p-2}
    \left(f_{n,1}^2+2f_{n,1}f_{n,2}\right)\right]\mathrm{d}x,
    \end{align}
    where the last inequality follows from the nonnegativity of
    \[
    |\nabla U|^{p-2} |\nabla \varphi_n|^2+(p-2)|\tilde{\omega}_n|^{p-2}f_{n,2}^2,
    \]
    thanks to \eqref{evbcb2i} and the fact that $f_{n,2}^2\leq |\nabla \varphi_n|^2$. Then, combining the convergence
    \begin{equation*}
    \begin{split}
    & \nabla \varphi_n \chi_{\mathcal{R}_n}\rightharpoonup 0,
    \quad f_{n,1}\to \frac{\nabla U}{|\nabla U|}\cdot\nabla \widehat{v},
    \quad f_{n,2}\chi_{\mathcal{R}_n}\rightharpoonup 0,
    \quad \mbox{locally in}\quad L^2(\mathbb{R}^N\backslash\{\mathbf{0}\}),\\
    & |\omega_n|\to |\nabla U|\quad \mbox{a.e.}, \quad |(B_R\backslash B_{\frac{1}{R}})\backslash\mathcal{R}_{n,R}|=|\mathcal{S}_{n,R}|\to 0,
    \end{split}
    \end{equation*}
    with the fact that
    \[
    |\tilde{\omega}_n|^{p-2}\leq C(p)|\nabla U|^{p-2},
    \]
    by Lebesgue's dominated convergence theorem, we deduce that
    \begin{align*}
    & \liminf_{n\to \infty}\int_{\mathcal{R}_{n,R}}
    |x|^{-\mu}\left[|\nabla U|^{p-2}\left(|\nabla \widehat{v}|^2+2 \nabla \varphi_n\cdot \nabla \widehat{v}\right)+(p-2)|\tilde{\omega}_n|^{p-2}
    \left(f_{n,1}^2+2f_{n,1}f_{n,2}\right)\right]\mathrm{d}x \\
    & \quad \to \int_{B_R\backslash B_{\frac{1}{R}}}
    |x|^{-\mu}\left[|\nabla U|^{p-2}|\nabla \widehat{v} |^2+(p-2)|\nabla U|^{p-2}\left(\frac{\nabla U\cdot\nabla \widehat{v} }{| \nabla U|}\right)^2\right]\mathrm{d}x,
    \end{align*}
    thus from \eqref{evbcbb2} we obtain
    \begin{align*}
    & \liminf_{n\to \infty}\int_{\mathcal{R}_{n,R}}
    |x|^{-\mu}\left[|\nabla U|^{p-2}|\nabla \widehat{v}_n|^2+(p-2)|\tilde{\omega}_n|^{p-2}\left(\frac{|\nabla (U+v_n)|-|\nabla U|}{\varepsilon_n}\right)^2\right]\mathrm{d}x \\
    & \quad \geq \int_{B_R\backslash B_{\frac{1}{R}}}
    |x|^{-\mu}\left[|\nabla U|^{p-2}|\nabla \widehat{v} |^2+(p-2)|\nabla U|^{p-2}\left(\frac{\nabla U\cdot\nabla \widehat{v} }{| \nabla U|}\right)^2\right]\mathrm{d}x.
    \end{align*}
    Recalling \eqref{evbcb2l} and \eqref{csirtc}, since $R>1$ is arbitrary and the integrand is nonnegative, this proves that
    \begin{align}\label{uinb2p2}
    & \int_{\mathbb{R}^N}
    |x|^{-\mu}\left[|\nabla U|^{p-2}|\nabla \widehat{v}|^2+(p-2)|\nabla U|^{p-2}\left(\frac{\nabla U\cdot\nabla \widehat{v}}{|\nabla U|}\right)^2\right]\mathrm{d}x \nonumber\\
    & \quad \leq \left[(r-1)+\tau\right]
    \int_{\mathbb{R}^N}|x|^{-s}U^{r-2}|\widehat{v}|^2\mathrm{d}x,
    \end{align}
    The orthogonality of $v_n$ (and also of $\widehat{v}_n$) implies that $\widehat{v}$ also is orthogonal to $T_{U} \mathcal{M}$. Since $\widehat{v}\in L^2_{s,*}(\mathbb{R}^N)$, \eqref{csirtl1} and \eqref{uinb2p2} contradict Proposition \ref{propevl}, the proof is complete.

    $\bullet$ {\em The case $\frac{2(N-\mu)}{N+2-s}< p<2$} which implies $r> 2$. If the statement fails, there exists a sequence $0\not\equiv v_n\to 0$ in $\mathcal{D}^{1,p}_{\mu}(\mathbb{R}^N)$, with $v_n$ orthogonal to $T_{U} \mathcal{M}$ in $L^2_{s,*}(\mathbb{R}^N)$, such that
    \begin{small}
    \begin{align}\label{evbc2g}
    & \int_{\mathbb{R}^N}
    |x|^{-\mu}\bigg[
    |\nabla U|^{p-2}|\nabla v_n|^2
    +(p-2)|\tilde{\omega}_n|^{p-2}(|\nabla (U+v_n)|-|\nabla U|)^2
    \nonumber\\ &\quad\quad +\gamma_0 \min\{|\nabla v_n|^p,|\nabla U|^{p-2}|\nabla v_n|^2\}
    \bigg]\mathrm{d}x \nonumber\\
    & \quad < \left[(r-1)+\tau\right]
    \int_{\mathbb{R}^N}|x|^{-s}U^{r-2}|v_n|^2\mathrm{d}x,
    \end{align}
    \end{small}
    where $\omega_n$ corresponds to $v_n$ as in the statement. As in the case $1<p<\frac{2(N-\mu)}{N+2-s}$, we define
    \[
    \varepsilon_n:=\left(\int_{\mathbb{R}^N}
    |x|^{-\mu}(|\nabla U|+|\nabla v_n|)^{p-2}|\nabla v_n|^2 \mathrm{d}x\right)^{\frac{1}{2}},\quad \widehat{v}_n=\frac{v_n}{\varepsilon_n},
    \]
    and we also have $\varepsilon_n\to 0$ as $n\to \infty$.

    Then, we also split $B_R\backslash B_{\frac{1}{R}}=\mathcal{R}_{n,R} \cup \mathcal{S}_{n,R}$, \eqref{evbcb2l} and \eqref{evbcb2li} hold also in this case, with the only difference that the last term in both equations now becomes
    \[
    \left[(r-1)+ \tau\right]\int_{\mathbb{R}^N}
    |x|^{-s}U^{r-2}|v_n|^2\mathrm{d}x,
    \]
    which is much simpler.

    We now observe that, by using H\"{o}lder inequality, we have
    \begin{small}\begin{align*}
    \int_{\mathbb{R}^N}
    |x|^{-\mu}|\nabla \widehat{v}_n|^p \mathrm{d}x
    \leq & \left(\int_{\mathbb{R}^N}
    |x|^{-\mu}(|\nabla U|+|\nabla v_n|)^{p-2}|\nabla \widehat{v}_n|^2 \mathrm{d}x\right)^{\frac{p}{2}}
    \left(\int_{\mathbb{R}^N}
    |x|^{-\mu}(|\nabla U|+|\nabla v_n|)^{p}\mathrm{d}x\right)^{1-\frac{p}{2}} \\
    = & \left(\int_{\mathbb{R}^N}
    |x|^{-\mu}(|\nabla U|+|\nabla v_n|)^{p}\mathrm{d}x\right)^{1-\frac{p}{2}} \\
    \leq & C(p)\left[\left(\int_{\mathbb{R}^N}
    |x|^{-\mu}|\nabla U|^{p}\mathrm{d}x\right)^{1-\frac{p}{2}}
    +\varepsilon_n^{\frac{p(2-p)}{2}}\left(\int_{\mathbb{R}^N} |x|^{-\mu}|\nabla \widehat{v}_n|^{p}\mathrm{d}x\right)^{1-\frac{p}{2}}\right]
    \end{align*}\end{small}
    from which it follows that
    \begin{equation}\label{vihl}
    \begin{split}
    \int_{\mathbb{R}^N}
    |x|^{-\mu}|\nabla \widehat{v}_n|^p \mathrm{d}x
    \leq C(N,p,\mu,s).
    \end{split}
    \end{equation}
    Thus, up to a subsequence, $\widehat{v}_n\rightharpoonup \widehat{v}$ in $\mathcal{D}^{1,p}_{\mu}(\mathbb{R}^N)$, then Rellich-Kondrachov compactness theorem implies $\widehat{v}_n \to \widehat{v}$ locally in $L^2(\mathbb{R}^N)$ (note that in this case $r>2$). In addition, Then by H\"{o}lder inequality and Sobolev inequality, together with \eqref{vihl}, yield for any $\rho>0$,
    \begin{align*}
    \int_{\mathbb{R}^N\backslash B_\rho} |x|^{-s}U^{r-2}|\widehat{v}_n|^2 \mathrm{d}x
    \leq & \left(\int_{\mathbb{R}^N\backslash B_\rho}|x|^{-s} U^{r}\mathrm{d}x\right)^{1-\frac{2}{r}}
    \left(\int_{\mathbb{R}^N\backslash B_\rho}|x|^{-s} |\widehat{v}_n|^{r}\mathrm{d}x\right)^{\frac{2}{r}} \\
    \leq & C\left(\int_{\mathbb{R}^N\backslash B_\rho}|x|^{-s} U^{r}\mathrm{d}x\right)^{1-\frac{2}{r}},
    \end{align*}
    while the last term tends to zero as $\rho\to +\infty$. Then by the strong convergence $\widehat{v}_n \to \widehat{v}$ locally in $L^2(\mathbb{R}^N)$, we conclude that $\widehat{v}_n\to \widehat{v}$ in $L^2_{s,*}(\mathbb{R}^N)$.

    In particular, letting $n\to \infty$ in the analogue of \eqref{evbcb2li} we obtain
    \begin{equation}\label{csirt0l}
    0<c(N,p,\mu,s,C_1,\gamma_0)
    \leq \|\widehat{v}\|_{L^2_{\beta,*}(\mathbb{R}^N)}
    \leq C(N,p,\mu,s,C_1).
    \end{equation}
    Similarly, the analogue of \eqref{evbcb2l} implies that
    \begin{equation}\label{csirt0g}
    |\mathcal{S}_{n,R}|\to 0\quad \mbox{and}\quad \int_{\mathbb{R}^N}
    |x|^{-\mu}|\nabla U|^{p-2}|\nabla \widehat{v}_n|^2 \mathrm{d}x\leq C(N,p,\mu,s),\quad \mbox{for all}\quad R>1.
    \end{equation}
    So, it follows from the weak convergence $\widehat{v}_n\rightharpoonup \widehat{v}$ in $\mathcal{D}^{1,p}_{\mu}(\mathbb{R}^N)$ that, up to a subsequence,
    \[
    \widehat{v}_n \chi_{\mathcal{R}_{n,R}}\rightharpoonup \widehat{v}\chi_{B_R\backslash B_{\frac{1}{R}}}
    \quad \mbox{locally in} \quad L^2(\mathbb{R}^N),\quad \mbox{for all}\quad R>1.
    \]
    Thanks to this bound, we can split
    \[
    \widehat{v}_n=\widehat{v}+\varphi_n,\quad\mbox{with}\quad \varphi_n:= \widehat{v}_n-\widehat{v},
    \]
    and very same argument as in the case $1<p<\frac{2(N-\mu)}{N+2-s}$ allows us to deduce that
    \begin{align*}
    & \liminf_{n\to \infty}\int_{\mathcal{R}_{n,R}}
    |x|^{-\mu}\left[|\nabla U|^{p-2}|\nabla \widehat{v}_n|^2+(p-2)|\tilde{\omega}_n|^{p-2}\left(\frac{|\nabla (U+v_n)|-|\nabla U|}{\varepsilon_n}\right)^2\right]\mathrm{d}x \\
    & \quad \geq \int_{B_R\backslash B_{\frac{1}{R}}}
    |x|^{-\mu}\left[|\nabla U|^{p-2}|\nabla \widehat{v} |^2+(p-2)|\nabla U|^{p-2}\left(\frac{\nabla U\cdot\nabla \widehat{v} }{| \nabla U|}\right)^2\right]\mathrm{d}x.
    \end{align*}
    Recalling \eqref{evbc2g}, since $R>1$ is arbitrary and the integrands above are nonnegative, this proves that \eqref{uinb2p2} holds, a contradiction to Proposition \ref{propevl} since $\widehat{v}$ is orthogonal to $T_{U} \mathcal{M}$ (being the strong $L^2_{s,*}(\mathbb{R}^N)$ limit of $\widehat{v}_n$).
    \end{proof}

\section{\bfseries Stability of (CKN) inequality}\label{sectpromr}

The main ingredient of the stability of (CKN) inequality \eqref{cknp} is contained in the two lemmas below, in which the behavior near the extremal functions set $\mathcal{M}=\{cU_{\lambda}: c\in\mathbb{R}, \lambda>0\}$ is studied.

    In order to shorten formulas, for each $u_n\in \mathcal{D}^{1,p}_{\mu}(\mathbb{R}^N)$ we denote
    \begin{equation*}
    \begin{split}
    \|u_n\|:
    =\left(\int_{\mathbb{R}^N}|x|^{-\mu}|\nabla u_n|^p \mathrm{d}x\right)^{\frac{1}{p}},
    \quad \|u_n\|_*: =\left(\int_{\mathbb{R}^N}|x|^{-s}|u_n|^{r} \mathrm{d}x\right)^{\frac{1}{r}},
    \end{split}
    \end{equation*}
    and
    \[
    d_n:=\inf_{v\in \mathcal{M}}
    \|u_n-v\|=\inf_{c\in\mathbb{R},\lambda>0}
    \|u_n-cU_{\lambda}\|.
    \]

    \begin{lemma}\label{lemma:rtnm2b}
    Assume that \eqref{cknib1c} holds and $p\geq 2$.
    There exists a small constant $\rho_1>0$ such that for any sequence $\{u_n\}\subset \mathcal{D}^{1,p}_{\mu}(\mathbb{R}^N)\backslash \mathcal{M}$ satisfying $\inf_n\|u_n\|>0$ and $d_n\to 0$,
    \begin{equation}\label{rtnmb}
    \liminf_{n\to\infty}
    \frac{\|u_n\|^p
    -\mathcal{S}\|u_n\|_*^p}
    {d_n^p}
    \geq \rho_1.
    \end{equation}
    \end{lemma}

    \begin{proof}
    We know that for each $u_n\in \mathcal{D}^{1,p}_{\mu}(\mathbb{R}^N)$, there exist $c_n\in\mathbb{R}$ and $\lambda_n>0$ such that $d_n=\|u_n-c_nU_{\lambda_n}\|$. In fact, since $p\geq 2$, for each fixed $n$, from Lemma \ref{lemui1p}, we obtain that for any $0<\kappa<1$, there exists a constant $\mathcal{C}_1=\mathcal{C}_1(p,\kappa)>0$ such that
    \begin{small}\begin{align}\label{ikeda}
    \|u_n-cU_\lambda\|^p
    = & \int_{\mathbb{R}^N}
    |x|^{-\mu}|\nabla u_n-c\nabla U_\lambda|^p \mathrm{d}x\nonumber\\
    \geq & \int_{\mathbb{R}^N}
    |x|^{-\mu}|\nabla u_n|^p \mathrm{d}x
    -pc\int_{\mathbb{R}^N}
    |x|^{-\mu}|\nabla u_n|^{p-2}  \nabla u_n\cdot \nabla U_\lambda \mathrm{d}x
    \nonumber\\ &
    +\mathcal{C}_1|c|^{p}\int_{\mathbb{R}^N}
    |x|^{-\mu}|\nabla U_\lambda|^{p} \mathrm{d}x
    +\frac{(1-\kappa)p}{2}c^2 \int_{\mathbb{R}^N}
    |x|^{-\mu}|\nabla u_n|^{p-2}  |\nabla U_\lambda|^2\mathrm{d}x \nonumber\\
    & +\frac{(1-\kappa)p(p-2)}{2}\int_{\mathbb{R}^N}
    |x|^{-\mu}|\bar{\omega}(\nabla u_n, \nabla u_n-c \nabla U_{\lambda})|^{p-2}(|\nabla u_n|-|\nabla u_n-c \nabla U_{\lambda}|)^2  \mathrm{d}x \nonumber\\
    \geq & \|u_n\|^p+ \mathcal{C}_1|c|^{p}\|U\|^p-pc\int_{\mathbb{R}^N}
    |x|^{-\mu}|\nabla u_n|^{p-2}  \nabla u_n\cdot \nabla U_\lambda \mathrm{d}x\nonumber\\
    \geq & \|u_n\|^p+ \mathcal{C}_1|c|^{p}\|U\|^p-p|c|\|U\|\|u_n\|^{p-1},
    \end{align}\end{small}
    where $\bar{\omega}:\mathbb{R}^{2N}\to \mathbb{R}^N$ corresponds to $\nabla u_n$ and $\nabla u_n-c \nabla U_{\lambda}$ as in Lemma \ref{lemui1p} for the case $p\geq 2$.
    Thus the minimizing sequence of $d_n$, say $\{c_{n,m},\lambda_{n,m}\}$, must satisfy $|c_{n,m}|\leq C$ which means $\{c_{n,m}\}$ is bounded.
    On the other hand,
    \begin{align*}
    \left|\int_{|\lambda x|\leq \rho}
    |x|^{-\mu}|\nabla u_n|^{p-2}  \nabla u_n\cdot \nabla U_\lambda \mathrm{d}x\right|
    \leq & \int_{|y|\leq \rho}
    |x|^{-\mu}|\nabla  (u_n)_{\frac{1}{\lambda}}(y)|^{p-1}|\nabla U(y)| \mathrm{d}y \\
    \leq & \|u_n\|^{p-1}\left(\int_{|y|\leq \rho}
    |x|^{-\mu}|\nabla U|^p \mathrm{d}y\right)^{\frac{1}{p}} \\
    = & o_\rho(1)
    \end{align*}
    as $\rho\to 0$ which is uniform for $\lambda>0$, where $(u_n)_{\frac{1}{\lambda}}(y)=\lambda^{-\frac{N-p-\mu}{p}}u_n(\lambda^{-1}y)$, and
    \begin{align*}
    \left|\int_{|\lambda x|\geq \rho}
    |x|^{-\mu}|\nabla u_n|^{p-2}  \nabla u_n\cdot \nabla U_\lambda \mathrm{d}x\right|
    \leq & \|U\|\left(\int_{|x|\geq \frac{\rho}{\lambda}}
    |x|^{-\mu}|\nabla u_n|^p \mathrm{d}y\right)^{\frac{1}{p}}
    =  o_\lambda(1)
    \end{align*}
    as $\lambda\to 0$ for any fixed $\rho>0$. By taking $\lambda\to 0$ and then $\rho\to 0$, we obtain
    \[\left|\int_{\mathbb{R}^N}
    |x|^{-\mu}|\nabla u_n|^{p-2}  \nabla u_n\cdot \nabla U_\lambda \mathrm{d}x\right| \to 0\quad \mbox{as}\quad \lambda\to 0.\]
    Moreover, by the explicit from of $U_\lambda$ we have
    \begin{equation*}
    \left|\int_{|\lambda x|\leq R}
    |x|^{-\mu}|\nabla u_n|^{p-2}  \nabla u_n\cdot \nabla U_\lambda \mathrm{d}x\right|
    \leq  \|U\|\left(\int_{| x|\leq \frac{R}{\lambda}}
    |x|^{-\mu}|\nabla u_n|^p \mathrm{d}x\right)^{\frac{1}{p}}
    = o_\lambda(1)
    \end{equation*}
    as $\lambda\to +\infty$ for any fixed $R>0$, and
    \begin{equation*}
    \begin{split}
    \left|\int_{|\lambda x|\geq R}
    |x|^{-\mu}|\nabla u_n|^{p-2}  \nabla u_n\cdot \nabla U_\lambda \mathrm{d}x\right|
    \leq & \int_{|y|\geq R}
    |x|^{-\mu}|\nabla (u_n)_{\frac{1}{\lambda}}(y)|^{p-1}|\nabla U(y)| \mathrm{d}y \\
    \leq & \|u_n\|^{p-1}\left(\int_{|y|\geq R}
    |x|^{-\mu}|\nabla U|^p \mathrm{d}y\right)^{\frac{1}{p}}
    =  o_R(1)
    \end{split}
    \end{equation*}
    as $R\to +\infty$ which is uniform for $\lambda>0$. Thus, by taking first $\lambda\to +\infty$ and then $R\to +\infty$, we also obtain
    \[
    \left|\int_{\mathbb{R}^N}
    |x|^{-\mu}|\nabla u_n|^{p-2}  \nabla u_n\cdot \nabla U_\lambda \mathrm{d}x\right| \to 0\quad \mbox{as}\quad \lambda\to +\infty.
    \]
    It follows from (\ref{ikeda}) and $d_n\to 0$, $\inf_n\|u_n\|>0$ that the minimizing sequence $\{c_{n,m},\lambda_{n,m}\}$ must satisfy $1/C\leq \lambda_{n,m}\leq C$ for some $C>1$ independent of $m$,  which means $\{\lambda_{n,m}\}$ is bounded. Thus for each $u_n\in  \mathcal{D}^{1,p}_{\mu}(\mathbb{R}^N)\backslash \mathcal{M}$, $d_n$ can be attained by some $c_n\in\mathbb{R}$ and $\lambda_n>0$.

    Since $\mathcal{M}$ is two-dimensional manifold embedded in $\mathcal{D}^{1,p}_{\mu}(\mathbb{R}^N)$, that is
    \[
    (c,\lambda)\in\mathbb{R}\times\mathbb{R}_+\to cU_\lambda\in \mathcal{D}^{1,p}_{\mu}(\mathbb{R}^N),
    \]
    then from \eqref{czkj}, under suitable transform, the tangential space is
    \[
    T_{c_n U_{\lambda_n}}\mathcal{M}={\rm Span}\left\{U_{\lambda_n}, \ \frac{\partial U_\lambda}{\partial\lambda}\Big|_{\lambda=\lambda_n}\right\}.
    \]
    Then we have,
    \[
    \int_{\mathbb{R}^N}|x|^{-\mu}|\nabla U_{\lambda_n}|^{p-2}\nabla \xi\cdot \nabla (u_n-c_n U_{\lambda_n}) \mathrm{d}x
    =0,\quad \mbox{for all}\quad \xi \in T_{c_n U_{\lambda_n}}\mathcal{M},
    \]
    particularly, taking $\xi=U_{\lambda_n}$ we obtain
    \begin{equation*}
    \int_{\mathbb{R}^N}|x|^{-\mu}|\nabla U_{\lambda_n}|^{p-2}\nabla U_{\lambda_n}\cdot \nabla (u_n-c_n U_{\lambda_n}) \mathrm{d}x
    =\int_{\mathbb{R}^N}|x|^{-s}U_{\lambda_n}^{r-1}(u_n-c_n U_{\lambda_n})\mathrm{d}x
    =0.
    \end{equation*}
    Let
    \begin{equation}\label{defunwn}
    u_n=c_n U_{\lambda_n}+d_n w_n,
    \end{equation}
     then $w_n$ is perpendicular to $T_{c_n U_{\lambda_n}}\mathcal{M}$, we have
    \begin{equation*}
    \|w_n\|=1,
    \end{equation*}
    and
    \begin{equation*}
    \int_{\mathbb{R}^N}|x|^{-s}U_{\lambda_n}^{r-1}w_n\mathrm{d}x
    =\int_{\mathbb{R}^N}|x|^{-\mu}|\nabla U_{\lambda_n}|^{p-2}\nabla U_{\lambda_n}\cdot \nabla w_n \mathrm{d}x=0.
    \end{equation*}
    From Lemma \ref{lemui1p}, for any $\kappa>0$, there exists a constant $\mathcal{C}_1=\mathcal{C}_1(p,\kappa)>0$ such that
    \begin{small}
    \begin{align}\label{epknug}
    \|u_n\|^p
    = & \int_{\mathbb{R}^N}
    |x|^{-\mu}|c_n \nabla U_{\lambda_n}+d_n \nabla w_n|^p \mathrm{d}x\nonumber\\
    \geq & |c_n|^{p}\int_{\mathbb{R}^N}
    |x|^{-\mu}|\nabla U_{\lambda_n}|^p \mathrm{d}x
    +p|c_n|^{p-2}c_nd_n \int_{\mathbb{R}^N}
    |x|^{-\mu}|\nabla U_{\lambda_n}|^{p-2}  \nabla U_{\lambda_n}\cdot \nabla w_n \mathrm{d}x  \nonumber\\
    & +\mathcal{C}_1d_n^{p}\int_{\mathbb{R}^N}
    |x|^{-\mu}|\nabla w_n|^{p} \mathrm{d}x
    +\frac{(1-\kappa)p}{2} |c_n|^{p-2}d_n^2\int_{\mathbb{R}^N}
    |x|^{-\mu}|\nabla U_{\lambda_n}|^{p-2}  |\nabla w_n|^2\mathrm{d}x \nonumber\\
    & +\frac{(1-\kappa)p(p-2)}{2}\int_{\mathbb{R}^N}
    |x|^{-\mu}|\bar{\omega}(c_n \nabla U_{\lambda_n},\nabla u_n)|^{p-2}(|c_n \nabla U_{\lambda_n}|-|\nabla u_n|)^2  \mathrm{d}x \nonumber\\
    = & |c_n|^{p}\|U\|^p+ \mathcal{C}_1d_n^{p}
     +\frac{(1-\kappa)p}{2} |c_n|^{p-2}d_n^2\int_{\mathbb{R}^N}
     |x|^{-\mu}|\nabla U_{\lambda_n}|^{p-2}  |\nabla w_n|^2\mathrm{d}x \nonumber\\
    & +\frac{(1-\kappa)p(p-2)}{2}\int_{\mathbb{R}^N}
    |x|^{-\mu}|\bar{\omega}(c_n \nabla U_{\lambda_n},\nabla u_n)|^{p-2}(|c_n \nabla U_{\lambda_n}|-|\nabla u_n|)^2  \mathrm{d}x.
    \end{align}
    \end{small}
    Then from Lemma \ref{lemui1p*l}, there exists a constant $\mathcal{C}_2=\mathcal{C}_2(r,\kappa)>0$ such that
    \begin{align*}
    \|u_n\|^{r}_*= & \int_{\mathbb{R}^N}|x|^{-s}|c_n U_{\lambda_n}+d_nw_n|^{r}  \mathrm{d}x\\
    \leq & |c_n|^{r}
    \int_{\mathbb{R}^N}|x|^{-s}U_{\lambda_n}^{r} \mathrm{d}x
    +|c_n|^{r-2}c_n r d_n \int_{\mathbb{R}^N}|x|^{-s}U_{\lambda_n}^{r-1}w_n \mathrm{d}x  \\
    & +\left(\frac{r(r-1)}{2}+\kappa\right)
    |c_n|^{r-2} d_n^2
    \int_{\mathbb{R}^N}|x|^{-s}U_{\lambda_n}^{r-2}w_n^2 \mathrm{d}x \\
    & +\mathcal{C}_2d_n^{r}
    \int_{\mathbb{R}^N}|x|^{-s}|w_n|^{r} \mathrm{d}x \\
    = & |c_n|^{r}\|U\|^p
    +\left(\frac{r(r-1)}{2}+\kappa\right)
    |c_n|^{r-2} d_n^2
    \int_{\mathbb{R}^N}|x|^{-s}
    U_{\lambda_n}^{r-2}w_n^2 \mathrm{d}x
    + o(d_n^p),
    \end{align*}
    since assumption \eqref{cknib1c} implies $p<r$.
    Thus, by the concavity of $t\mapsto t^{\frac{p}{r}}$, we have
    \begin{align}\label{epkeyiyxbb}
    \|u_n\|^p_*
    \leq &  |c_n|^p\|U\|^{\frac{p^2}{r}}
    + \frac{p|c_n|^{p-2} d_n^2}{r}
    \left(\frac{r(r-1)}{2}+\kappa\right)
    \|U\|^{\frac{p^2}{r-p}}
    \int_{\mathbb{R}^N}|x|^{-s}U_{\lambda_n}^{r-2}w_n^2
    \mathrm{d}x
    +o(d_n^p).
    \end{align}
    Note also that Lemma \ref{lemsgap} implies
    \begin{align}\label{epknugb}
    \int_{\mathbb{R}^N}|x|^{\alpha} \left[|\nabla c_nU_{\lambda_n}|^{p-2}  |\nabla d_nw_n|^2
    +(p-2)|\bar{\omega}(c_n \nabla U_{\lambda_n},\nabla u_n)|^{p-2}(|c_n \nabla U_{\lambda_n}|-|\nabla u_n|)^2\right]  dx& \nonumber \\
    \geq \left[(p^*_{\alpha,\beta}-1)+\tau\right]|c_n|^{p-2}
    \int_{\mathbb{R}^N}|x|^{\beta}
    U_{\lambda_n}^{p^*_{\alpha,\beta}-2}|d_nw_n|^2dx& ,
    \end{align}
    for some constant $\tau>0$.
    Therefore, as $d_n\to 0$, combining \eqref{epknug} with \eqref{epkeyiyxbb}, it follows from \eqref{epknugb} that, by choosing $\kappa>0$ small enough,
    \begin{align*}
    \|u_n\|^p
    -\mathcal{S}\|u_n\|_*^p
    \geq & |c_n|^{p}\|U\|^p+ \mathcal{C}_1d_n^{p}
    +\frac{(1-\kappa)p}{2} d_n^2\int_{\mathbb{R}^N}
    |x|^{-\mu}|\nabla c_n U_{\lambda_n}|^{p-2}  |\nabla w_n|^2\mathrm{d}x \\
    & +\frac{(1-\kappa)p(p-2)}{2}\int_{\mathbb{R}^N}
    |x|^{-\mu}|\bar{\omega}(c_n \nabla U_{\lambda_n},\nabla u_n)|^{p-2}(|c_n \nabla U_{\lambda_n}|-|\nabla u_n|)^2 \mathrm{d}x \\
    & -\mathcal{S}\Bigg\{|c_n|^p\|U\|^{\frac{p^2}{r}} +o(d_n^p) \\
    & \quad \quad+ \frac{p|c_n|^{r-2} d_n^2}{r}
    \left(\frac{r(r-1)}{2}+\kappa\right)
    \|U\|^{\frac{p^2}{r}-p}
    \int_{\mathbb{R}^N}|x|^{-s}U_{\lambda_n}^{r-2}w_n^2 \mathrm{d}x\Bigg\} \\
    \geq  & \mathcal{C}_1d_n^p -o(d_n^p)
    + \frac{(1-\kappa)p  |c_n|^{p-2} d_n^2}{2}\left[(r-1)+\tau\right]
    \int_{\mathbb{R}^N}|x|^{-s}U_{\lambda_n}^{r-2}w_n^2 \mathrm{d}x\\
    & - \frac{p|c_n|^{p-2} d_n^2}{r}
    \left(\frac{r(r-1)}{2}+\kappa\right)
    \mathcal{S}\|U\|^{\frac{p^2}{r}-p}
    \int_{\mathbb{R}^N}|x|^{-s}U_{\lambda_n}^{r-2}w_n^2 \mathrm{d}x \\
    \geq  & \mathcal{C}_1d_n^p -o(d_n^p),
    \end{align*}
    due to $\|U\|^p=\mathcal{S}^{\frac{r}{r-p}}$, then (\ref{rtnmb}) follows immediately.
    \end{proof}

    \begin{lemma}\label{lemma:rtnm2b2}
    Assume that \eqref{cknib1c} holds and $1<p<2$.
    There exists a small constant $\rho_2>0$ such that for any sequence $\{u_n\}\subset \mathcal{D}^{1,p}_{\mu}(\mathbb{R}^N)\backslash \mathcal{M}$ satisfying $\inf_n\|u_n\|>0$ and $d_n\to 0$,
    \begin{equation}\label{rtnmb2}
    \liminf_{n\to\infty}
    \frac{\|u_n\|^p
    -\mathcal{S}\|u_n\|_*^p}
    {d_n^2}
    \geq \rho_2.
    \end{equation}
    \end{lemma}

    \begin{proof}
    We know that for each $u_n\in \mathcal{D}^{1,p}_{\mu}(\mathbb{R}^N)$, there exist $c_n\in\mathbb{R}$ and $\lambda_n>0$ such that $d_n=\|u_n-c_nU_{\lambda_n}\|$.
    In fact, since $1< p<2$, for each fixed $n$, from Lemma \ref{lemui1p}, we obtain that for any $0<\kappa<1$, there exists a constant $\mathcal{C}_1=\mathcal{C}_1(p,\kappa)>0$ such that
    \begin{small}
    \begin{align}\label{ikeda2}
    \|u_n-cU_\lambda\|^p
    = & \int_{\mathbb{R}^N}
    |x|^{-\mu}|\nabla u_n-c\nabla U_\lambda|^p \mathrm{d}x\nonumber\\
    \geq & \int_{\mathbb{R}^N}
    |x|^{-\mu}|\nabla u_n|^p \mathrm{d}x
    -pc\int_{\mathbb{R}^N}
    |x|^{-\mu}|\nabla u_n|^{p-2}  \nabla u_n\cdot \nabla U_\lambda \mathrm{d}x  \nonumber\\ &
    +\mathcal{C}_1|c|^{2}\int_{\mathbb{R}^N}
    |x|^{-\mu}\min\{|c|^{p-2}|\nabla U_\lambda|^p,|\nabla u_n|^{p-2}|\nabla U_\lambda|^2\} \mathrm{d}x
    \nonumber\\ &
    +\frac{(1-\kappa)p}{2}c^2 \int_{\mathbb{R}^N}
    |x|^{-\mu}|\nabla u_n|^{p-2}  |\nabla U_\lambda|^2\mathrm{d}x \nonumber\\
    & +\frac{(1-\kappa)p(p-2)}{2}\int_{\mathbb{R}^N}
    |x|^{-\mu}|\tilde{\omega}(\nabla u_n, \nabla u_n-c \nabla U_{\lambda})|^{p-2}(|\nabla u_n|-|\nabla u_n-c \nabla U_{\lambda}|)^2  \mathrm{d}x \nonumber\\
    \geq & \int_{\mathbb{R}^N}
    |x|^{-\mu}|\nabla u_n|^p \mathrm{d}x
    -pc\int_{\mathbb{R}^N}
    |x|^{-\mu}|\nabla u_n|^{p-2}  \nabla u_n\cdot \nabla U_\lambda \mathrm{d}x  \nonumber\\ &
    +\mathcal{C}_1|c|^{2}\int_{\mathbb{R}^N}
    |x|^{-\mu}\min\{|c|^{p-2}|\nabla U_\lambda|^p,|\nabla u_n|^{p-2}|\nabla U_\lambda|^2\} \mathrm{d}x
    \nonumber\\
    \geq & \|u_n\|^p-p|c|\|U\|\|u_n\|^{p-1}\nonumber\\ &
    +\mathcal{C}_1|c|^{2}\int_{\mathbb{R}^N}
    |x|^{-\mu}\min\{|c|^{p-2}|\nabla U_\lambda|^p,|\nabla u_n|^{p-2}|\nabla U_\lambda|^2\} \mathrm{d}x,
    \end{align}
    \end{small}
    where $\tilde{\omega}:\mathbb{R}^{2N}\to \mathbb{R}^N$ corresponds to $\nabla u_n$ and $\nabla u_n-c \nabla U_{\lambda}$ the same as $\tilde{\omega}$ in Lemma \ref{lemui1p} for the case $1<p<2$, since from \eqref{evbcb2i} it holds that
    \begin{align*}
    0\leq & c^2 \int_{\mathbb{R}^N}
    |x|^{-\mu}|\nabla u_n|^{p-2}  |\nabla U_\lambda|^2\mathrm{d}x
    \nonumber\\ & + (p-2)\int_{\mathbb{R}^N}
    |x|^{-\mu}|\tilde{\omega}(\nabla u_n, \nabla u_n-c \nabla U_{\lambda})|^{p-2}(|\nabla u_n|-|\nabla u_n-c \nabla U_{\lambda}|)^2  \mathrm{d}x.
    \end{align*}
    Therefore the minimizing sequence of $d_n$, say $\{c_{n,m},\lambda_{n,m}\}$, must satisfy $|c_{n,m}|\leq C$ for some $C>0$ independent of $m$,  which means $\{c_{n,m}\}$ is bounded.
    On the other hand, taking the same steps as those in Lemma \ref{lemma:rtnm2b}, we deduce that
    \[
    \left|\int_{\mathbb{R}^N}
    |x|^{-\mu}|\nabla u_n|^{p-2}  \nabla u_n\cdot \nabla U_\lambda \mathrm{d}x\right|
    \to 0\quad \mbox{as}\quad \lambda\to 0,
    \]
    and
    \[
    \left|\int_{\mathbb{R}^N}
    |x|^{-\mu}|\nabla u_n|^{p-2}  \nabla u_n\cdot \nabla U_\lambda \mathrm{d}x\right|
    \to 0\quad \mbox{as}\quad \lambda\to +\infty.
    \]
    It follows from (\ref{ikeda}) and $d_n\to 0$, $\inf_n\|u_n\|>0$ that the minimizing sequence $\{c_{n,m},\lambda_{n,m}\}$ must satisfy $1/C\leq |\lambda_{n,m}|\leq C$ for some $C>1$ independent of $m$,  which means $\{\lambda_{n,m}\}$ is bounded. Thus for each $u_n\in  \mathcal{D}^{1,p}_{\mu}(\mathbb{R}^N)\backslash \mathcal{M}$, $d_n$ can also be attained by some $c_n\in\mathbb{R}$ and $\lambda_n>0$.

    As stated in Lemma \ref{lemma:rtnm2b}, we have
    \[
    T_{c_n U_{\lambda_n}}\mathcal{M}={\rm Span}\left\{U_{\lambda_n}, \ \frac{\partial U_\lambda}{\partial\lambda}\Big|_{\lambda=\lambda_n}\right\},
    \]
    and
    \begin{equation*}
    \int_{\mathbb{R}^N}|x|^{-s}U_{\lambda_n}^{p-1}(u_n-c_n U_{\lambda_n})\mathrm{d}x
    =\int_{\mathbb{R}^N}
    |x|^{-\mu}|\nabla U_{\lambda_n}|^{p-2}\nabla U_{\lambda_n}\cdot \nabla (u_n-c_n U_{\lambda_n}) \mathrm{d}x
    =0.
    \end{equation*}
    Let
    \begin{equation}\label{defunwn2}
    u_n=c_n U_{\lambda_n}+d_n w_n,
    \end{equation}
    then $w_n$ is perpendicular to $T_{c_n U_{\lambda_n}}\mathcal{M}$, we have
    \begin{equation*}
    \|w_n\|=1,
    \end{equation*}
    and
    \begin{equation*}
    \int_{\mathbb{R}^N}|x|^{-s}U_{\lambda_n}^{r-1}w_n\mathrm{d}x
    =\int_{\mathbb{R}^N}|x|^{-\mu}|\nabla U_{\lambda_n}|^{p-2}\nabla U_{\lambda_n}\cdot \nabla w_n \mathrm{d}x=0.
    \end{equation*}
    Since $1<p<2$, from Lemma \ref{lemui1p} there exists a constant $\mathcal{C}_1=\mathcal{C}_1(p,\kappa)>0$ such that
    \begin{align}\label{epknug2}
    \|u_n\|^p
    = & \int_{\mathbb{R}^N}
    |x|^{-\mu}|c_n \nabla U_{\lambda_n}+d_n \nabla w_n|^p \mathrm{d}x \nonumber\\
    \geq & |c_n|^{p}\int_{\mathbb{R}^N}
    |x|^{-\mu}|\nabla U_{\lambda_n}|^p \mathrm{d}x
    +p|c_n|^{p-2}c_nd_n \int_{\mathbb{R}^N}
    |x|^{-\mu}|\nabla U_{\lambda_n}|^{p-2}  \nabla U_{\lambda_n}\cdot \nabla w_n \mathrm{d}x  \nonumber\\
    & +\frac{(1-\kappa)p}{2} |c_n|^{p-2}d_n^2\int_{\mathbb{R}^N}
    |x|^{-\mu}|\nabla U_{\lambda_n}|^{p-2}  |\nabla w_n|^2\mathrm{d}x \nonumber\\
    & +\frac{(1-\kappa)p(p-2)}{2}\int_{\mathbb{R}^N}
    |x|^{-\mu}|\tilde{\omega}(c_n \nabla U_{\lambda_n},\nabla u_n)|^{p-2}(|c_n \nabla U_{\lambda_n}|-|\nabla u_n|)^2  \mathrm{d}x \nonumber\\
    & +\mathcal{C}_1d_n^{2} \int_{\mathbb{R}^N}
    |x|^{-\mu}\min\{d_n^{p-2}|\nabla w_n|^{p}, |c_n \nabla U_{\lambda_n}|^{p-2}|\nabla w_n|^{2}\} \mathrm{d}x  \nonumber\\
    = & |c_n|^{p}\|U\|^p
    +\frac{(1-\kappa)p}{2} |c_n|^{p-2}d_n^2\int_{\mathbb{R}^N}
    |x|^{-\mu}|\nabla U_{\lambda_n}|^{p-2}  |\nabla w_n|^2\mathrm{d}x \nonumber\\
    & +\frac{(1-\kappa)p(p-2)}{2}\int_{\mathbb{R}^N}
    |x|^{-\mu}|\tilde{\omega}(c_n \nabla U_{\lambda_n},\nabla u_n)|^{p-2}(|c_n \nabla U_{\lambda_n}|-|\nabla u_n|)^2  \mathrm{d}x\nonumber\\
    & +\mathcal{C}_1d_n^{2} \int_{\mathbb{R}^N}
    |x|^{-\mu}\min\{d_n^{p-2}|\nabla w_n|^{p}, |c_n \nabla U_{\lambda_n}|^{p-2}|\nabla w_n|^{2}\} \mathrm{d}x.
    \end{align}
    Then we consider the following two cases:

    $\bullet$ {\em The case $1<p\leq\frac{2(N-\mu)}{N+2-s}$} which implies $r\leq2$.

    From Lemma \ref{lemui1p*l}, for any $\kappa>0$, there exists a constant $\mathcal{C}_2=\mathcal{C}_2(r,\kappa)>0$ such that
    \begin{align*}
    \|u_n\|^{r}_*
    = & \int_{\mathbb{R}^N}|x|^{-s}|c_n U_{\lambda_n}+d_nw_n|^{r}  \mathrm{d}x\\
    \leq & |c_n|^{r}
    \int_{\mathbb{R}^N}|x|^{-s}U_{\lambda_n}^{r} \mathrm{d}x
    +r |c_n|^{r-2}c_n d_n \int_{\mathbb{R}^N}|x|^{-s}U_{\lambda_n}^{r-1}w_n \mathrm{d}x  \\
    & +\left(\frac{r(r-1)}{2}+\kappa\right)d_n^2
    \int_{\mathbb{R}^N}|x|^{-s}\frac{(|c_n U_{\lambda_n}|+\mathcal{C}_2|d_nw_n|)^{r}}{|c_n U_{\lambda_n}|^2+|d_nw_n|^2}w_n^2 \mathrm{d}x \\
    = & |c_n|^{r}\|U\|^p
    +\left(\frac{r(r-1)}{2}+\kappa\right)d_n^2
    \int_{\mathbb{R}^N}|x|^{-s}\frac{(|c_n U_{\lambda_n}|+\mathcal{C}_2|d_nw_n|)^{r}}{|c_n U_{\lambda_n}|^2+|d_nw_n|^2}w_n^2 \mathrm{d}x.
    \end{align*}
    Thus, by the concavity of $t\mapsto t^{\frac{p}{r}}$, we have
    \begin{align}\label{epkeyiyxbb2}
    \|u_n\|^p_*
    \leq  &  |c_n|^p\|U\|^{\frac{p^2}{r}}
    \nonumber\\ & +\frac{p}{r}
    \left(\frac{r(r-1)}{2}+\kappa\right)d_n^2
    \|U\|^{\frac{p^2}{r}-p}
    \int_{\mathbb{R}^N}|x|^{-s}\frac{(|c_n U_{\lambda_n}|+\mathcal{C}_2|d_nw_n|)^{r}}{|c_n U_{\lambda_n}|^2+|d_nw_n|^2}w_n^2 \mathrm{d}x.
    \end{align}
    Therefore, as $d_n\to 0$, combining \eqref{epkeyiyxbb2} with \eqref{epknug2}, it follows from Lemma \ref{lemsgap} that, by choosing $\kappa>0$ small enough,
    \begin{small}
    \begin{align*}
    \|u_n\|^p
    -\mathcal{S}\|u_n\|_*^p
    \geq & |c_n|^{p}\|U\|^p
    +\frac{(1-\kappa)p}{2}d_n^2\int_{\mathbb{R}^N}
    |x|^{-\mu}|\nabla c_n U_{\lambda_n}|^{p-2}  |\nabla w_n|^2\mathrm{d}x \\
    & +\frac{(1-\kappa)p(p-2)}{2}\int_{\mathbb{R}^N}
    |x|^{-\mu}|\tilde{\omega}(c_n \nabla U_{\lambda_n},\nabla u_n)|^{p-2}(|c_n \nabla U_{\lambda_n}|-|\nabla u_n|)^2 \mathrm{d}x\\
    & +\mathcal{C}_1d_n^{2} \int_{\mathbb{R}^N}
    |x|^{-\mu}\min\{d_n^{p-2}|\nabla w_n|^{p}, |c_n \nabla U_{\lambda_n}|^{p-2}|\nabla w_n|^{2}\}
    \mathrm{d}x \\
    & -\mathcal{S}\left[|c_n|^p\|U\|^{\frac{p^2}{r}}
    +\left(\frac{p(r-1)}{2}
    +\frac{p\kappa}{r}\right)d_n^2
    \|U\|^{\frac{p^2}{r}-p}
    \int_{\mathbb{R}^N}|x|^{-s}\frac{(|c_n U_{\lambda_n}|+\mathcal{C}_2|d_nw_n|)^{r}}{|c_n U_{\lambda_n}|^2+|d_nw_n|^2}w_n^2 \mathrm{d}x\right] \\
    =  & \frac{(1-\kappa)p}{2} d_n^2\int_{\mathbb{R}^N}
    |x|^{-\mu}|\nabla c_nU_{\lambda_n}|^{p-2}  |\nabla w_n|^2\mathrm{d}x
    \\ & +\frac{(1-\kappa)p(p-2)}{2}\int_{\mathbb{R}^N}
    |x|^{-\mu}|\tilde{\omega}(c_n \nabla U_{\lambda_n},\nabla u_n)|^{p-2}(|c_n \nabla U_{\lambda_n}|-|\nabla u_n|)^2  \mathrm{d}x\\
    & +\mathcal{C}_1d_n^{2} \int_{\mathbb{R}^N}
    |x|^{-\mu}\min\{d_n^{p-2}|\nabla w_n|^{p}, |c_n \nabla U_{\lambda_n}|^{p-2}|\nabla w_n|^{2}\} \mathrm{d}x
    \\ & -\left(\frac{p(r-1)}{2}
    +\frac{p\kappa}{r}\right)d_n^2
    \int_{\mathbb{R}^N}|x|^{-s}\frac{(|c_n U_{\lambda_n}|+\mathcal{C}_2|d_nw_n|)^{r}}{|c_n U_{\lambda_n}|^2+|d_nw_n|^2}w_n^2 \mathrm{d}x,
    \end{align*}
    \end{small}
    since $\|U\|^p=\mathcal{S}^{\frac{r}{r-p}}$. Lemma \ref{lemsgap2} allows us to reabsorb the last term above: more precisely, we have
    \begin{small}
    \begin{align*}
    & \|u_n\|^p
    -\mathcal{S}\|u_n\|_*^p \\
    \geq & p\left(\frac{(1-\kappa)}{2} - \frac{(r-1)+\frac{2}{r}\kappa}
    {2(r-1)
    +2\tau}\right) \\
    & \quad \times
    \int_{\mathbb{R}^N}
    |x|^{-\mu}\left[|\nabla c_nU_{\lambda_n}|^{p-2}  |\nabla d_nw_n|^2
    +(p-2)|\tilde{\omega}(c_n \nabla U_{\lambda_n},\nabla u_n)|^{p-2}\left(|c_n \nabla U_{\lambda_n}|-|\nabla u_n|\right)^2\right]\mathrm{d}x \\
    & + d_n^{2}\left(\mathcal{C}_1 - \gamma_0\frac{p\left[(r-1)
    +\frac{2}{r}\kappa\right]}
    {2(r-1)+2\tau}\right)
    \int_{\mathbb{R}^N}
    |x|^{-\mu}\min\{d_n^{p-2}|\nabla w_n|^{p}, |c_n \nabla U_{\lambda_n}|^{p-2}|\nabla w_n|^{2}\} \mathrm{d}x.
    \end{align*}
    \end{small}
    Now, let us recall the definition of $\omega$, as stated below Lemma \ref{lemui1p}, we have
    \[
    |\nabla c_nU_{\lambda_n}|^{p-2}  |\nabla d_nw_n|^2
    +(p-2)|\tilde{\omega}(c_n \nabla U_{\lambda_n},\nabla u_n)|^{p-2}\left(|c_n \nabla U_{\lambda_n}|-|\nabla u_n|\right)^2\geq 0,
    \]
    then choosing $\kappa>0$ small enough such that
    \[
    \frac{(1-\kappa)}{2} - \frac{(r-1)+\frac{2}{r}\kappa}{2(r-1)
    +2\tau }\geq 0,
    \]
    and then $\gamma_0>0$ small enough such that
    \[
    \frac{\mathcal{C}_1}{2}\geq \gamma_0\frac{p\left[(r-1)
    +\frac{2}{r}\kappa\right]}
    {2(r-1)+2\tau},
    \]
    we eventually arrive at
    \begin{align}\label{emn21}
    & \|u_n\|^p
    -\mathcal{S}\|u_n\|_*^p
    \geq \frac{\mathcal{C}_1}{2}d_n^{2}\int_{\mathbb{R}^N}
    |x|^{-\mu}\min\{d_n^{p-2}|\nabla w_n|^{p}, |c_n \nabla U_{\lambda_n}|^{p-2}|\nabla w_n|^{2}\} \mathrm{d}x.
    \end{align}
    Observe that, since $1<p<2$, it follows by H\"{o}lder inequality that
    \begin{small}\begin{align*}
    \left(\int_{\{d_n|\nabla w_n|\geq|c_n \nabla U_{\lambda_n}|\}}
    |x|^{-\mu}|\nabla w_n|^{p} \mathrm{d}x\right)^{\frac{2}{p}}
    \leq & \left(\int_{\{d_n|\nabla w_n|\geq|c_n \nabla U_{\lambda_n}|\}}
    |x|^{-\mu}|\nabla U_{\lambda_n}|^{p} \mathrm{d}x\right)^{\frac{2}{p}-1}  \\
    & \quad \times
    \int_{\{d_n|\nabla w_n|\geq|c_n \nabla U_{\lambda_n}|\}}
    |x|^{-\mu}|\nabla U_{\lambda_n}|^{p-2}|\nabla w_n|^{2} \mathrm{d}x  \\
    \leq & \mathcal{S}^{\frac{r(\frac{2}{p}-1)}
    {r-p}} \int_{\{d_n|\nabla w_n|\geq|c_n \nabla U_{\lambda_n}|\}}
    |x|^{-\mu}|\nabla U_{\lambda_n}|^{p-2}|\nabla w_n|^{2}
    \mathrm{d}x,
    \end{align*}\end{small}
    then we obtain
    \begin{small}\begin{align}\label{emn21b}
    & \int_{\mathbb{R}^N}
    |x|^{-\mu}\min\{d_n^{p-2}|\nabla w_n|^{p}, |c_n \nabla U_{\lambda_n}|^{p-2}|\nabla w_n|^{2}\} \mathrm{d}x
    \nonumber\\
    & \quad = d_n^{p-2}\int_{\{d_n|\nabla w_n|<|c_n \nabla U_{\lambda_n}|\}}
    |x|^{-\mu}|\nabla w_n|^{p} \mathrm{d}x
    +\int_{\{d_n|\nabla w_n|\geq |c_n \nabla U_{\lambda_n}|\}}
    |x|^{-\mu}|c_n \nabla U_{\lambda_n}|^{p-2}|\nabla w_n|^{2}\mathrm{d}x
    \nonumber\\
    & \quad \geq d_n^{p-2}\int_{\{d_n|\nabla w_n|< |c_n \nabla U_{\lambda_n}|\}}
    |x|^{-\mu}|\nabla w_n|^{p} \mathrm{d}x
    + c\left(\int_{\{d_n|\nabla w_n|\geq|c_n \nabla U_{\lambda_n}|\}}
    |x|^{-\mu}|\nabla w_n|^{p} \mathrm{d}x\right)^{\frac{2}{p}}
    \nonumber\\
    & \quad \geq c\left(\int_{\mathbb{R}^N}
    |x|^{-\mu}|\nabla w_n|^{p} \mathrm{d}x\right)^{\frac{2}{p}}=c,
    \end{align}\end{small}
    for some constant $c>0$.
    The conclusion (\ref{rtnmb2}) follows immediately from \eqref{emn21} and \eqref{emn21b}.

    $\bullet$ {\em The case $\frac{2(N-\mu)}{N+2-s}< p<2$} which implies $r>2$.

    The proof is very similar to the previous case, with very small changes. From Lemma \ref{lemui1p*l}, we have that for any $\kappa>0$, there exists a constant $\mathcal{C}_2=\mathcal{C}_2(r,\kappa)>0$ such that
    \begin{align*}
    \|u_n\|_*^{r}
    = & \int_{\mathbb{R}^N}|x|^{-s}|u_n|^{r} \mathrm{d}x \\
    \leq & |c_n|^{r}
    \int_{\mathbb{R}^N}|x|^{-s}U_{\lambda_n}^{r} \mathrm{d}x
    +|c_n|^{r-2}c_n r d_n \int_{\mathbb{R}^N}|x|^{-s}U_{\lambda_n}^{r-1}w_n \mathrm{d}x  \\
    & +\left(\frac{r(r-1)}{2}+\kappa\right)
    |c_n|^{r-2} d_n^2
    \int_{\mathbb{R}^N}|x|^{-s}U_{\lambda_n}^{r-2}w_n^2
    +\mathcal{C}_2d_n^{r}
    \int_{\mathbb{R}^N}|x|^{-s}|w_n|^{r} \mathrm{d}x \\
    = & |c_n|^{r}\|U\|^p
    +\left(\frac{r(r-1)}{2}+\kappa\right)
    |c_n|^{r-2} d_n^2
    \int_{\mathbb{R}^N}|x|^{-s}
    U_{\lambda_n}^{r-2}w_n^2 \mathrm{d}x
    + o(d_n^2).
    \end{align*}
    Then by the concavity of $t\mapsto t^{\frac{p}{r}}$, we have
    \begin{align}\label{epkeyiyxbbg2}
    \|u_n\|_*^{p}
    \leq  &  |c_n|^p\|U\|^{\frac{p^2}{r}}
    + o(d_n^2)
    + \frac{p|c_n|^{p-2} d_n^2}{r}
    \left(\frac{r(r-1)}{2}+\kappa\right)
    \|U\|^{\frac{p^2}{r}-p}
    \int_{\mathbb{R}^N}|x|^{-s}U_{\lambda_n}^{r-2}w_n^2 \mathrm{d}x.
    \end{align}
    Hence, arguing as in the case $1<p<\frac{2(N-\mu)}{N+2-s}$,
    Therefore, as $d_n\to 0$, combining \eqref{epknug2} with \eqref{epkeyiyxbbg2}, it follows from Lemma \ref{lemsgap} that, by choosing $\kappa>0$ small enough,
    \begin{small}
    \begin{align*}
    \|u_n\|^p
    -\mathcal{S}\|u_n\|_*^p
    \geq & |c_n|^{p}\|U\|^p
    +\frac{(1-\kappa)p}{2} |c_n|^{p-2}d_n^2\int_{\mathbb{R}^N}
    |x|^{-\mu}|\nabla U_{\lambda_n}|^{p-2}  |\nabla w_n|^2\mathrm{d}x \\
    & +\frac{(1-\kappa)p(p-2)}{2}\int_{\mathbb{R}^N}
    |x|^{-\mu}|\tilde{\omega}(c_n \nabla U_{\lambda_n},\nabla u_n)|^{p-2}(|c_n \nabla U_{\lambda_n}|-|\nabla u_n|)^2  \mathrm{d}x\\
    & +\mathcal{C}_1d_n^{2} \int_{\mathbb{R}^N}
    |x|^{-\mu}\min\{d_n^{p-2}|\nabla w_n|^{p}, |c_n \nabla U_{\lambda_n}|^{p-2}|\nabla w_n|^{2}\} \mathrm{d}x \\
    & -\mathcal{S}\Bigg\{|c_n|^p\|U\|^{\frac{p^2}{r}}
    + o(d_n^2) \\
    & \quad \quad+ \frac{p|c_n|^{p-2} d_n^2}{r}
    \left(\frac{r(r-1)}{2}+\kappa\right)
    \|U\|^{\frac{p^2}{r}-p}
    \int_{\mathbb{R}^N}|x|^{-s}U_{\lambda_n}^{r-2}
    w_n^2 \mathrm{d}x
    \Bigg\} \\
    \geq  & \frac{(1-\kappa)p}{2} |c_n|^{p-2}d_n^2\int_{\mathbb{R}^N}
    |x|^{-\mu}|\nabla U_{\lambda_n}|^{p-2}  |\nabla w_n|^2\mathrm{d}x \\
    & +\frac{(1-\kappa)p(p-2)}{2}\int_{\mathbb{R}^N}
    |x|^{-\mu}|\tilde{\omega}(c_n \nabla U_{\lambda_n},\nabla u_n)|^{p-2}(|c_n \nabla U_{\lambda_n}|-|\nabla u_n|)^2  \mathrm{d}x\\
    & +\mathcal{C}_1d_n^{2} \int_{\mathbb{R}^N}
    |x|^{-\mu}\min\{d_n^{p-2}|\nabla w_n|^{p}, |c_n \nabla U_{\lambda_n}|^{p-2}|\nabla w_n|^{2}\} \mathrm{d}x \\
    & -\left(\frac{p(r-1)}{2}
    +\frac{p\kappa}{r}\right)d_n^2
    \int_{\mathbb{R}^N}|x|^{-s}U_{\lambda_n}^{r-2}w_n^2 \mathrm{d}x
    - o(d_n^2).
    \end{align*}
    \end{small}
    Lemma \ref{lemsgap2} allows us to reabsorb the last term above: more precisely, we have
    \begin{small}
    \begin{align*}
    & \|u_n\|^p
    -\mathcal{S}\|u_n\|_*^p
    \\ \geq & p\left(\frac{(1-\kappa)}{2} - \frac{(r-1)+\frac{2}{r}\kappa}{2(r-1)
    +2\tau }\right) \\
    & \quad \times
    \int_{\mathbb{R}^N}
    |x|^{-\mu}\left[|\nabla U_{\lambda_n}|^{p-2}  |\nabla d_nw_n|^2
    +(p-2)|\tilde{\omega}(c_n \nabla U_{\lambda_n},\nabla u_n)|^{p-2}\left(|c_n \nabla U_{\lambda_n}|-|\nabla u_n|\right)^2\right]\mathrm{d}x \\
    & + d_n^{2}\left(\mathcal{C}_1 - \gamma_0\frac{p\left[(r-1)
    +\frac{2}{r}\kappa\right]}
    {2(r-1)+2\tau}\right)
    \int_{\mathbb{R}^N}
    |x|^{-\mu}\min\{d_n^{p-2}|\nabla w_n|^{p}, |c_n \nabla U_{\lambda_n}|^{p-2}|\nabla w_n|^{2}\} \mathrm{d}x
    - o(d_n^2).
    \end{align*}
    \end{small}
    Note that $\tau>0$ is a constant, we can always choose $\kappa>0$ small enough such that
    \[
    \frac{(1-\kappa)}{2} - \frac{(r-1)+\frac{2}{r}\kappa}{2(r-1)
    +2\tau}\geq 0,
    \]
    and $\gamma_0>0$ small enough such that
    \[
    \frac{\mathcal{C}_1}{2}\geq \gamma_0\frac{p\left[(r-1)+\frac{2}{r}\kappa\right]}
    {2(r-1)+2\tau}.
    \]
    From \eqref{emn21b}, we eventually arrive at
    \begin{small}
    \begin{align*}
    \|u_n\|^p
    -\mathcal{S}\|u_n\|_*^p
    \geq & \frac{\mathcal{C}_1}{2}d_n^{2}\int_{\mathbb{R}^N}
    |x|^{-\mu}\min\{d_n^{p-2}|\nabla w_n|^{p}, |c_n \nabla U_{\lambda_n}|^{p-2}|\nabla w_n|^{2}\} \mathrm{d}x
    - o(d_n^2)
    \geq cd_n^{2},
    \end{align*}
    \end{small}
    for some constant $c>0$, thus the conclusion (\ref{rtnmb2}) follows immediately.
    \end{proof}

Now, we are ready to prove our main result.
\subsection{\bf Proof of Theorem \ref{thmprtp}}\label{subspmr}

    By homogeneity, we can assume that $\|u\|_*=1$. Since $\inf_{v\in \mathcal{M}}\|u-v\|\leq \|u\|$ which implies $
\inf_{v\in \mathcal{M}}\left(\frac{\|u-v\|}{\|u\|}\right)\leq 1$,
it is suffices to prove the result of $\|u\|-\mathcal{S}^{\frac{1}{p}} \ll 1$.

  Now, we argue by contradiction. In fact, if the theorem is false then there exists a sequence $\{u_n\}\subset \mathcal{D}^{1,p}_{\mu}(\mathbb{R}^N)\backslash \mathcal{M}$ satisfying $\|u_n\|_*=1$ and $\|u_n\|-\mathcal{S}^{\frac{1}{p}}
    \ll 1$ such that
    \begin{align*}
    \liminf_{n\to\infty}
    \frac{\|u_n\|-\mathcal{S}^{\frac{1}{p}}}
    {d_n^\gamma\|u_n\|^{-\gamma}}
    \to 0,\quad \mbox{as}\quad n\to \infty,
    \end{align*}
    where $d_n=\inf_{v\in \mathcal{M}}\|u_n-v\|$ and $\gamma=\max\{2,p\}$. Observe that, the assumption $ \|u_n\|-\mathcal{S}^{\frac{1}{p}} \ll 1$ combines with $\|u_n\|\geq \mathcal{S}^{\frac{1}{p}}\|u_n\|_*=\mathcal{S}^{\frac{1}{p}}$ indicate $\{\|u_n\|\}$ is bounded away from zero and infinity and $\|u_n\|-\mathcal{S}^{\frac{1}{p}} \geq c(\|u_n\|^p
    -\mathcal{S})$ for some constant $c>0$ thus
    \begin{equation*}
    \liminf_{n\to\infty}
    \frac{\|u_n\|^p
    -\mathcal{S}\|u_n\|_*^p}
    {d_n^\gamma}
    \to 0,\quad \mbox{as}\quad n\to \infty.
    \end{equation*}
    After selecting a subsequence we can assume that $d_n\to \varpi\in[0,\infty)$ since $d_n=\inf_{c\in\mathbb{R}, \lambda>0}\|u_n-cU_{\lambda}\|\leq \|u_n\|$. If $\varpi=0$, then we deduce a contradiction by Lemmas  \ref{lemma:rtnm2b} and \ref{lemma:rtnm2b2}.

    The other possibility only is that $\varpi>0$, that is
    $d_n\to \varpi>0$,
    then we must have
    \begin{equation}\label{wbsi}
    \|u_n\|^p
    -\mathcal{S}\to 0,\quad \|u_n\|_*=1.
    \end{equation}
    From the concentration and compactness principle established by Tan and Yang \cite{TY04}, we know that going if necessary to a subsequence, there is a sequence $\{\lambda_n\}\subset\mathbb{R}^+$ such that
    \begin{equation*}
    \lambda_n^{\frac{N-p-\mu}{p}}u_n(\lambda_n x)\to U_*\quad \mbox{in}\quad \mathcal{D}^{1,p}_\mu(\mathbb{R}^N)\quad \mbox{as}\quad n\to \infty,
    \end{equation*}
    for some $U_*\in\mathcal{M}$, which implies
    \begin{equation*}
    d_n=\mathrm{dist}(u_n,\mathcal{M})=\mathrm{dist}
    \left(\lambda_n^{\frac{N-p-\mu}{p}}u_n(\lambda_n x),\mathcal{M}\right)\to 0 \quad \mbox{as}\quad n\to \infty,
    \end{equation*}
    which leads to a contradiction.

    Therefore, the proof of Theorem \ref{thmprtp} is completed.
    \qed

\medskip

\noindent{\bfseries Acknowledgements}

The research bas been supported by National Natural Science Foundation of China (No. 12371121).

\appendix

\section{\bfseries Several crucial algebra inequalities}\label{sectpls}

\begin{lemma}\label{lemui1p}
    \cite[Lemmas 2.1]{FZ22} Let $x, y\in\mathbb{R}^N$. Then for any $\kappa>0$, there exists a constant $\mathcal{C}_1=\mathcal{C}_1(p,\kappa)>0$ such that
     the following inequalities hold.

$\bullet$ For $p\geq 2$,
    \begin{align}\label{uinb2p}
    |x+y|^p
    \geq & |x|^p+ p|x|^{p-2}x\cdot y+ \frac{1-\kappa}{2}\left(p|x|^{p-2}|y|^2+ p(p-2)|\bar{\omega}|^{p-2}(|x|-|x+y|)^2 \right) \nonumber\\
    & +\mathcal{C}_1 |y|^p ,
    \end{align}
    where
    \begin{eqnarray*}
    \bar{\omega}=\bar{\omega}(x,x+y)=
    \left\{ \arraycolsep=1.5pt
       \begin{array}{ll}
        x,\ \ &{\rm if}\ \ |x|<|x+y|,\\[3mm]
        \left(\frac{|x+y|}{|x|}\right)^{\frac{1}{p-2}}(x+y),\ \ &{\rm if}\ \  |x+y|\leq |x|.
        \end{array}
    \right.
    \end{eqnarray*}

    $\bullet$ For $1<p<2$,
    \begin{align}\label{uinb1p}
    |x+y|^p
    \geq & |x|^p+ p|x|^{p-2}x\cdot y+ \frac{1-\kappa}{2}\left[p|x|^{p-2}|y|^2+ p(p-2)|\tilde{\omega}|^{p-2}(|x|-|x+y|)^2 \right] \nonumber\\
    & +\mathcal{C}_1\min\{|y|^p,|x|^{p-2}|y|^2\},
    \end{align}
    where
    \begin{eqnarray*}
    \tilde{\omega}=\tilde{\omega}(x,x+y)=
    \left\{ \arraycolsep=1.5pt
       \begin{array}{ll}
        \left(\frac{|x+y|}{(2-p)|x+y|+(p-1)|x|}\right)^{\frac{1}{p-2}}x,\ \ &{\rm if}\ \ |x|<|x+y|,\\[3mm]
        x,\ \ &{\rm if}\ \  |x+y|\leq |x|.
        \end{array}
    \right.
    \end{eqnarray*}
    \end{lemma}

    Note that for all $y\in\mathbb{R}^N$, $|\bar{\omega}(x,x+y)|^{p-2}\leq |x|^{p-2}$ if $p\geq 2$, and
    $ |\tilde{\omega}(x,x+y)|\leq \frac{|x|^{p-2}}{2-p}$ if $1<p<2$. Furthermore, $|x|^{p-2}|y|^2+(p-2)|\tilde{\omega}|^{p-2}(|x|-|x+y|)^2\geq 0$, see \cite[(2.2)]{FZ22} for details.

    \begin{lemma}\label{lemui1p*l}
    Let $a, b\in\mathbb{R}$. Then for any $\kappa>0$, there exists a constant $\mathcal{C}_2=\mathcal{C}_2(r,\kappa)>0$ where $r=\frac{p(N-s)}{N-p-\mu}$ such that
     the following inequalities hold.

    $\bullet$ For $1<p\leq\frac{2(N-\mu)}{N+2-s}$,
    \begin{equation}\label{uinx2pl}
    |a+b|^{r}
    \leq |a|^{r}+ r|a|^{r-2}a b
    + \left(\frac{r(r-1)}{2}+\kappa\right)
    \frac{(|a|+\mathcal{C}_2|b|)^{r}}
    {|a|^2+|b|^2}|b|^{2}.
    \end{equation}

    $\bullet$ For $\frac{2(N-\mu)}{N+2-s}< p<N$,
    \begin{equation}\label{uinx2pb}
    |a+b|^{r}
    \leq |a|^{r}+ r|a|^{r-2}a b
    + \left(\frac{r(r-1)}{2}+\kappa\right)
    |a|^{{r}-2}|b|^2 +\mathcal{C}_2|b|^{r}.
    \end{equation}
    \end{lemma}

    \begin{proof}
    Those results can be directly followed from \cite[Lemma 3.2]{FN19} and \cite[Lemma 2.4]{FZ22}.
    \end{proof}


\begin{thebibliography}{99}

\bibitem{ACP05}
Abdellaoui, B., Colorado, E., Peral, I.: {\em Some improved Caffarelli-Kohn-Nirenberg inequalities}. Calc. Var. Partial Differential Equations {\bf 23}(3), 327--345 (2005)

\bibitem{Au76}
Aubin, T.: {\em Probl\`{e}mes isop\'{e}rimtriqu\'{e}s et espaces de Sobolev}. J. Differ. Geom. {\bf 11}, 573--598 (1976)

\bibitem{BWW03}
Bartsch, T., Weth T., Willem, M.: {\em A Sobolev inequality with remainder term and critical equations on domains with topology for the polyharmonic operator}. Calc. Var. Partial Dif. {\bf 18}, 253--268 (2003)

\bibitem{BE91}
Bianchi, G., Egnell, H.: {\em A note on the Sobolev inequality}. J. Funct. Anal. {\bf 100}(1), 18--24 (1991)

\bibitem{BCG21}
Bonheure, D., Casteras, J., Gladiali, F.: {\em Bifurcation analysis of the Hardy-Sobolev equation}. J. Diff. Equ. {\bf 296}, 759--798 (2021)

\bibitem{BrL85}
Brezis H., Lieb, E.H.: {\em Sobolev inequalities with remainder terms}. J. Funct. Anal. {\bf 62}, 73--86 (1985)

\bibitem{BW02}
Byeon, J., Wang, Z.-Q.: {\em Symmetry breaking of extremal functions for the Caffarelli-Kohn-Nirenberg inequalities}. Commun. Contemp. Math. {\bf 4}(3), 457--465 (2002)

\bibitem{CKN84}
Caffarelli, L., Kohn, R., Nirenberg, L.: {\em First order interpolation inequalities with weights}. Compos. Math. {\bf 53}, 259--275 (1984)

\bibitem{CM13}
Caldiroli, P., Musina, R.: {\em
Symmetry breaking of extremals for the Caffarelli-Kohn-Nirenberg inequalities in a non-Hilbertian setting}. Milan J. Math. {\bf 81}(2), 421--430 (2013)

\bibitem{CW01}
Catrina, F., Wang, Z.-Q.: {\em On the Caffarelli-Kohn-Nirenberg inequalities: sharp constants, existence (and nonexistence), and symmetry of extremal functions}. Comm. Pure Appl. Math. {\bf 54}(2), 229--258 (2001)

\bibitem{CFW13}
Chen, S., Frank, R.L., Weth, T.: {\em Remainder terms in the fractional Sobolev inequality}. Indiana Univ. Math. J. {\bf 62}(4), 1381--1397 (2013)

\bibitem{CFMP09}
Cianchi, A., Fusco, N., Maggi, F., Pratelli, A.: {\em The sharp Sobolev inequality in quantitative form}. J. Eur. Math. Soc. (JEMS) {\bf 11}(5), 1105--1139 (2009)

\bibitem{CC22}
Ciraolo, G., Corso, R.: {\em Symmetry for positive critical points of Caffarelli-Kohn-Nirenberg inequalities}. Nonlinear Anal. {\bf 216}, Paper No. 112683, 23 pp (2022)

\bibitem{Do21}
Dolbeault, J.: {\em Functional inequalities: nonlinear flows and entropy methods as a tool for obtaining sharp and constructive results}. Milan J. Math. {\bf 89}(2), 355--386 (2021)

\bibitem{DEL16}
Dolbeault, J., Esteban, M.J., Loss, M.: {\em Rigidity versus symmetry breaking via nonlinear flows on cylinders and Euclidean spaces}. Invent. Math. {\bf 206}(2), 397--440 (2016)

\bibitem{DELT09}
Dolbeault, J., Esteban, M.J., Loss, M., Tarantello, G.: {\em
On the symmetry of extremals for the Caffarelli-Kohn-Nirenberg inequalities}. Adv. Nonlinear Stud. {\bf 9}(4), 713--726 (2009)

\bibitem{DT23-1}
Deng, S., Tian, X.: {\em Some weighted fourth-order Hardy-H\'{e}non equations}. J. Funct. Anal. {\bf 284}(1), Paper No. 109745 (2023)

\bibitem{DT22}
Deng, S., Tian, X.: {\em Caffarelli-Kohn-Nirenberg-type inequalities related to weighted $p$-Laplace equations}. Preprint. \url{https://arxiv.org/abs/2212.05459} (2022)

\bibitem{DT23}
Deng, S., Tian, X.: {\em Stability of Hardy-Sobolev inequality involving $p$-Laplace}. Preprint. \url{https://arxiv.org/abs/2301.07442} (2023)

\bibitem{DT23R2}
Deng, S., Tian, X.: {\em
On the stability of Caffarelli-Kohn-Nirenberg inequality in $\mathbb{R}^2$}. Preprint. \url{https://arxiv.org/abs/2308.04111} (2023)

\bibitem{FS03}
Felli, V., Schneider, M.: {\em Perturbation results of critical elliptic equations of Caffarelli-Kohn-Nirenberg type}. J. Diff. Equ. {\bf 191}, 121--142 (2003)

\bibitem{FN19}
Figalli, A., Neumayer, R.: {\em Gradient stability for the Sobolev inequality: the case $p\geq 2$}. J. Eur. Math. Soc. (JEMS) {\bf 21}(2), 319--354 (2019)

\bibitem{FZ22}
Figalli, A., Zhang Y.: {\em Sharp gradient stability for the Sobolev inequality}. Duke Math. J. {\bf 171}(12), 2407--2459 (2022)

\bibitem{GY00}
Ghoussoub, N., Yuan, C.: {\em Multiple solutions for quasi-linear PDEs involving the critical Sobolev and Hardy exponents}. Trans. Amer. Math. Soc. {\bf 352}(12), 5703--5743 (2000)

\bibitem{Ho97}
Horiuchi, T.: {\em Best constant in weighted Sobolev inequality with weights being powers of distance from the origin}. J. Inequal. Appl. {\bf 1}(3), 275--292 (1997)

\bibitem{LL17}
Lam, N., Lu, G.: {\em Sharp constants and optimizers for a class of Caffarelli-Kohn-Nirenberg inequalities}. Adv. Nonlinear Stud. {\bf 17}(3), 457--480 (2017)

\bibitem{Li85-2}
Lions, P.L.: {\em The concentration-compactness principle in the calculus of variations. The limit case. \uppercase\expandafter{\romannumeral 2}}.  Rev. Mat. Iberam. {\bf 1}(2), 45--121 (1985)

\bibitem{LW00}
Lu, G., Wei, J.: {\em  On a Sobolev inequality with remainder terms}, Proc. Amer. Math. Soc. {\bf 128}(1), 75--84 (2000).

\bibitem{Ma85}
Maz'ja, V. G.: {\em Sobolev spaces}. Translated from the Russian by T. O. Shaposhnikova. Springer Series in Soviet Mathematics. Springer-Verlag, Berlin, 1985.

\bibitem{Mu14}
Musina, R.: {\em Weighted Sobolev spaces of radially symmetric functions}. Ann. Mat. Pura Appl. (4) {\bf 193}(6), 1629--1659 (2014)

\bibitem{Ne20}
Neumayer, R.: {\em A note on strong-form stability for the Sobolev inequality}. Calc. Var. Partial Differential Equations {\bf 59}(1), Paper No. 25 (2020)

\bibitem{Ng19}
Nguyen, V.H.: {\em The sharp Gagliardo-Nirenberg-Sobolev inequality in quantitative form}. J. Funct. Anal. {\bf 277}(7), 2179--2208 (2019)

\bibitem{PV21}
Pistoia, A., Vaira, G.: {\em Nondegeneracy of the bubble for the critical p-Laplace equation}. Proc. Roy. Soc. Edinburgh Sect. A {\bf 151}(1), 151--168 (2021)

\bibitem{RSW02}
R\u{a}dulescu, V., Smets, D., Willem, M.: {\em Hardy-Sobolev inequalities with remainder terms}. Topol. Methods Nonlinear Anal. {\bf 20}(1), 145--149 (2002)

\bibitem{ST18}
Sano, M., Takahashi, F.: {\em Some improvements for a class of the Caffarelli-Kohn-Nirenberg inequalities}. Differ. Integral Equ. {\bf 31}(1-2), 57--74 (2018)

\bibitem{SW03}
Smets, D., Willem, M.: {\em Partial symmetry and asymptotic behavior for some elliptic variational problems}. Calc. Var. Partial Differential Equations {\bf 18}(1), 57--75 (2003)

\bibitem{SW22}
Su, J., Wang, C.:{\em Weighted critical exponents of Sobolev-type embeddings for radial functions}. Adv. Nonlinear Stud. {\bf 22}(1), 143--158 (2022)

\bibitem{Ta76}
Talenti, G.: {\em Best constant in Sobolev inequality}. Ann. Mat. Pura Appl. {\bf 110}, 353--372 (1976)

\bibitem{TY04}
Tan, J., Yang, J.: {\em On the singular variational problems}. Acta Math. Sci. Ser. B (Engl. Ed.) {\bf 24}(4), 672--690 (2004)

\bibitem{WaWi00}
Wang, Z.-Q., Willem, M.: {\em Singular minimization problems}. J. Differential Equations {\bf 161}(2), 307--320 (2000)

\bibitem{WaWi03}
Wang, Z.-Q., Willem, M.: {\em Caffarelli-Kohn-Nirenberg inequalities with remainder terms}. J. Funct. Anal. {\bf 203}(2), 550--568 (2003)

\bibitem{WW22}
Wei, J., Wu, Y.: {\em On the stability of the Caffarelli-Kohn-Nirenberg inequality}. Math. Ann. {\bf 384}(3-4), 1509--1546 (2022)

\bibitem{ZZ15}
Zhong, X., Zou, W.: {\em Existence of extremal functions for a family of Caffarelli-Kohn-Nirenberg inequalities}. Preprint. \url{https://arxiv.org/abs/1504.00433} (2015)

\end{thebibliography}
    \end{document}